\documentclass[11pt,reqno]{amsart}
\usepackage[a4paper,margin=25mm]{geometry}
\usepackage[T1]{fontenc}
\usepackage{lmodern}
\usepackage{amsmath,amssymb,amsthm,mathtools,bm,microtype,enumitem}
\usepackage[hidelinks,hypertexnames=false]{hyperref}
\hypersetup{pdftitle={Quadratic splitting for a weighted Monge-Ampere equation on a product domain},pdfauthor={Jiahuan Li, Yilu Liu, Xi-Nan Ma, Dekai Zhang}}
\allowdisplaybreaks
\setlist{itemsep=2pt,topsep=4pt}
\newtheorem{theorem}{Theorem}[section]
\newtheorem{proposition}[theorem]{Proposition}
\newtheorem{lemma}[theorem]{Lemma}
\newtheorem{corollary}[theorem]{Corollary}
\theoremstyle{definition}

\theoremstyle{remark}\newtheorem{remark}[theorem]{Remark}
\numberwithin{equation}{section}
\newcommand{\R}{\mathbb R}

\newcommand{\cl}{\overline}
\newcommand{\loc}{\mathrm{loc}}

\newcommand{\ee}{\bm e}

\DeclareMathOperator{\tr}{tr}

\DeclareMathOperator{\dist}{dist}

\newlist{proofsteps}{enumerate}{1}
\setlist[proofsteps]{label=\textbf{Step \arabic*.},ref=\arabic*,leftmargin=0pt,itemindent=3.6em,labelindent=0pt,labelwidth=3.1em,labelsep=.5em,listparindent=\parindent,itemsep=\medskipamount,parsep=0pt,topsep=\medskipamount}
\title[Quadratic splitting on a product domain]{Quadratic splitting for a weighted Monge--Amp\`ere equation on $(0,\infty)^l\times\R^m$}
\author{Jiahuan Li}
\author{Yilu Liu}
\author{Xi-Nan Ma}
\author{Dekai Zhang}
\date{}

\begin{document}
\begin{abstract}
We study a Monge--Amp\`ere equation arising from the existence
problem for complete K\"ahler--Einstein metrics with zero Ricci
curvature on divisor complements. We prove a splitting theorem for
nonnegative smooth strictly convex solutions on a product of positive
half lines and a Euclidean space. The potential and its first
derivatives extend continuously to the boundary. The directional
derivative in the equation is positive in the interior and vanishes
on the boundary. After an affine change of variables, the solution
is the sum of a function on the product of half lines and a positive
definite quadratic form in the Euclidean variables. No growth
condition or homogeneity is assumed.
\end{abstract}
\maketitle
\tableofcontents

\section{Introduction}\label{intro:chapter}
The J\"orgens--Calabi--Pogorelov theorem
\cite{Jorgens1954,Calabi1958,Pogorelov1972} states that an entire
classical convex solution of $\det D^2u=1$ is a quadratic polynomial.
Jian and Wang \cite{JianWang2013} and Savin
\cite{Savin2013Boundary,Savin2014} studied classification and
boundary regularity in a half space.
For transport between half spaces with densities equal to powers of
the normal coordinate, Jhaveri and Savin \cite{JS}
proved a classification theorem.

In this paper we consider the equation in \cite{CTY} on
$\Omega=(0,\infty)^l\times\R^m$, where $l,m\ge1$. We prove that the solution is quadratic in $y$
after an affine change of variables. We do not assume a growth
condition or homogeneity of the solution.
Let $\ee=(1,\ldots,1)^T\in\R^l$. The notation
$C^1_{\loc}(\cl\Omega)$ requires the function and its first
derivatives to extend continuously to the boundary, locally on
bounded sets. Our main result is as follows.
\begin{theorem}\label{rf:splitting}

Let $k>0$ and let
$\Phi\in C^\infty(\Omega)\cap C^1_{\loc}(\cl\Omega)$ satisfy
$\Phi\ge0$ and $D^2\Phi>0$ in $\Omega$.
Set $W=\sum_{i=1}^l\Phi_{x_i}$ and assume
\begin{equation}\label{intro:real-equation}
 \begin{gathered}
 W>0,\\
 W^k\det D^2\Phi=1\quad \hbox{in }\Omega,
 \\
 W=0\quad \hbox{on }\partial\Omega.
 \end{gathered}
\end{equation}
There are constant matrices $A=A^T>0$, $B\in\R^{m\times l}$, and a constant vector $b\in\R^m$, with $B\ee=0$, such that
\begin{equation}\label{intro:split-equation}
 \Phi(x,y)=\psi(x)+\frac{1}{2}(y+Bx+b)^TA(y+Bx+b).
\end{equation}
Here $\psi=\inf_y\Phi\ge0$ is strictly convex, belongs to
$C^\infty((0,\infty)^l)\cap C^1_{\loc}([0,\infty)^l)$, and satisfies
\begin{equation}\label{intro:residual-equation}
 \begin{gathered}
 \Bigl(\sum_{i=1}^l\psi_{x_i}\Bigr)^k\det D^2\psi
 =\frac{1}{\det A},\\
 \sum_{i=1}^l\psi_{x_i}>0
 \quad \hbox{in }(0,\infty)^l,
 \end{gathered}
\end{equation}
with $\sum_i\psi_{x_i}=0$ on the entire boundary. Conversely, any choice of $\psi,A,B,b$ with these properties defines a solution of \eqref{intro:real-equation} by \eqref{intro:split-equation}.

\end{theorem}

Differentiating \eqref{intro:split-equation}, we have
$W=\sum_i\psi_{x_i}+(B\ee)^TA(y+Bx+b)$ and $D_yW=AB\ee$.
Since $W>0$ for all $y\in\R^m$, we have $B\ee=0$.
For $l=2$ and $m=1$, the shear need not vanish: if $\psi$
satisfies \eqref{intro:residual-equation} with $A=1$, then for any
$a\in\R$ the potential
\[
 \Phi(x_1,x_2,y)=\psi(x_1,x_2)
                 +\frac{1}{2}\bigl(y+a(x_1-x_2)\bigr)^2
\]
has the form \eqref{intro:split-equation}. Thus the affine term $Bx$
cannot in general be omitted. The classification of $\psi$ in
\eqref{intro:residual-equation} remains to be studied.

\subsection{Rigidity on a half space and related results}
Jian and Wang \cite{JianWang2013} studied the Dirichlet problem
for $\det D^2u=\bigl(\frac{u_n}{x_n}\bigr)^{n+2}$ in $\{x_n>0\}$
with quadratic boundary values. Jian and Wang
\cite{JianXianduoWang2023} proved a Liouville theorem for
viscosity solutions under upper and lower quadratic growth bounds,
and considered perturbations of the equation and boundary values.
For convex solutions of $\det D^2u=x_n^\alpha$ in $\{x_n>0\}$
with $u(x',0)=\frac{|x'|^2}{2}$, Savin \cite{Savin2014} proved
that, if $\alpha\ge0$ and $u=O(|x|^{3+\alpha-\varepsilon})$ for
some $\varepsilon>0$, then
\[
 u(Sx)=b x_n+\frac{1}{2}|x'|^2
       +\frac{x_n^{2+\alpha}}{(1+\alpha)(2+\alpha)}
\]
for a sliding $S$ preserving $x_n=0$. Savin and Zhang
\cite{SavinZhang2020} studied the singular density $x_n^{-a}$,
$0<a<1$, under upper and lower anisotropic growth bounds.
Jia \cite{Jia2025Degenerate} considered compact perturbations of a
density equal to a positive power of the normal coordinate.

For the Neumann problem $\det D^2u=1$ in $\R^n_+$ with
$u_n=ax_1$ on $\{x_n=0\}$, Jian and Tu \cite{JianTu2023} proved that the solution
is quadratic when $n=2$. For $n>2$, the same conclusion holds for
zero Neumann values or under quadratic growth on a subspace of
dimension $n-2$. Wang and Zhou \cite{WangZhou2024} obtained
classification results in dimension two for related degenerate or
singular Dirichlet and Neumann problems.

Huang, Tang and Wang \cite{HuangTangWang2024Obstacle} proved a
Bernstein theorem for a singular equation in a half space. For the
operator with normal entry $\psi_{nn}+\frac{b\psi_n}{x_n}$, $b>-1$, they
obtained quadratic rigidity under zero normal derivative, global
$C^{1,1}$ control and uniform ellipticity.

Jhaveri and Savin \cite{JS} classified potentials transporting
measures with densities equal to powers of the normal coordinate
between half spaces, with the gradient image prescribed. For $l=1$, the two density exponents
in our equation are $0$ and $k$. We derive the gradient image from
the equation in Proposition~\ref{extR:globalpair}. For $l>1$, the boundary of $\Omega$ has intersecting faces. The condition $D_{\ee}\Phi=0$ is
imposed on all faces, and is not a separate normal derivative
condition on each face.

Boundary regularity for the Dirichlet problem was studied by Ivochkina
\cite{Ivochkina1980Boundary}, Krylov \cite{Krylov1983Domain},
Caffarelli, Nirenberg and Spruck
\cite{CaffarelliNirenbergSpruck1984Dirichlet,CaffarelliNirenbergSpruck1987Correction},
and Trudinger and Wang \cite{TrudingerWang2008Boundary}.
Interior estimates and the totally degenerate Dirichlet problem
were studied in
\cite{TrudingerUrbas1984SecondDerivative,CaffarelliNirenbergSpruck1986Degenerate}.
Savin \cite{Savin2012BoundaryLocalization,Savin2013Boundary}
proved boundary localization and pointwise estimates.
Hong, Huang and Wang \cite{HongHuangWang2011}, Savin
\cite{Savin2014}, and Le and Savin \cite{LeSavin2017} studied
boundary regularity for densities vanishing at the boundary.

Boundary estimates for Neumann and oblique problems on bounded
domains were proved in
\cite{LionsTrudingerUrbas1986Neumann,Urbas1998Oblique,JiangTrudingerXiang2016Neumann,ShiJiang2021Neumann,ShiJiang2021AugmentedNeumann}.
Jian and Tu \cite{JianTu2026Oblique} proved boundary
$C^{2,\alpha}$ estimates for the oblique Monge--Amp\`ere problem.

Caffarelli and Li \cite{CaffarelliLi2003} studied quadratic
asymptotics for exterior problems in the whole space. In a half space,
Jia, Li and Li \cite{JiaLiLi2020} proved quadratic asymptotics under
upper and lower quadratic growth bounds and quadratic values on
the flat boundary. Boundary perturbations were considered in
\cite{JiaLi2021,JiaLiXuemei2023}, and density perturbations with
prescribed asymptotics in \cite{Jia2023HalfSpace}. Ma and Jia
\cite{MaJia2024} assumed the growth bound
$O(|x|^{3-\varepsilon})$, $\varepsilon>0$. Liang
\cite{Liang2026HalfSpaceExpansion} obtained an expansion of higher
order under upper and lower quadratic growth bounds. Li and Liu
\cite{LiLiu2026HalfSpace} studied exterior asymptotics for fully
nonlinear elliptic equations in half spaces.

Le and Savin \cite{LeSavin2023Polygon}
classified convex solutions of the Dirichlet problem with constant
density on the planar first quadrant. Their classification includes nonquadratic solutions for
some densities. They also proved polygonal boundary estimates
assuming a $C^2$ convex strict subsolution. Huang and Shen
\cite{HuangShen2025Cones} studied expansions and higher regularity
at conical points. For transport on squares and hypercubes,
Jhaveri \cite{Jhaveri2019Instability} and Chen, Liu and Wang
\cite{ChenLiuWang2023Hypercubes} used reflection and, in dimension
two, a partial Legendre transform.

Collins and Firester \cite{CollinsFirester2025} constructed
homogeneous transport potentials from a cone to a half space with
density $s^k$. Collins, Firester and Tong
\cite{CollinsFiresterTong2025} constructed homogeneous transports
between cones under obliqueness and density assumptions.
For potentials homogeneous of degree two transporting Lebesgue
measure between cones, they proved quadratic splitting along
translation subspaces by Brascamp--Lieb rigidity
\cite{CollinsFiresterTong2026}. Here we consider the density $s^k$ without assuming homogeneity
of the potential.

Equation~\eqref{intro:real-equation} also arises in the construction
of complete K\"ahler metrics with vanishing Ricci curvature on
divisor complements. Earlier constructions of such metrics are
due to Calabi \cite{Calabi1979} and Tian and Yau
\cite{TianYau1990I,TianYau1991II}. When the divisor classes are proportional, the
generalized Calabi ansatz of Collins and Li \cite{CollinsLi2024}
and Collins, Tong and Yau \cite{CTY} leads to
\eqref{intro:real-equation} in logarithmic normal coordinates.
At a relative interior point of a face of codimension $l$ of a
simple polytope, the closed tangent cone is linearly isomorphic to
$[0,\infty)^l\times\R^m$, with the Euclidean factor tangent to the face \cite{Schneider2014}.

\subsection{The proof}
We first prove that the gradient image is a half space and
establish the corresponding mass identities. We use the section
estimates of Caffarelli
\cite{Caffarelli1990Localization,Caffarelli1991Regularity,Caffarelli1992Boundary,Caffarelli1996BoundaryII}
and the arguments of Caffarelli and McCann \cite{CM} and Jhaveri
and Savin \cite{JS}. We verify the doubling condition for the
Alexandrov measure and apply the engulfing lemma in
\cite{CM}.
For further estimates on sections,
see \cite{CaffarelliGutierrez1996RealAnalysis,GutierrezHuang2000Sections,LeNguyen2013BoundarySections}.

For the normalized potentials in Lemma~\ref{extR:polarity}, we prove
$(f_H)_{yy}\le C(1+f_H)$ and $(v_H)_{\eta\eta}\le C(1+v_H)$
for $\eta\perp e_s$, $|\eta|=1$, with $C$ independent of $H$. Tangential estimates for related
equations were proved in \cite{JS} and
\cite{CT}. Here we combine interior Pogorelov,
Evans--Krylov and Schauder estimates
\cite{DPF,Evans1982,Krylov1983,GilbargTrudinger2001} with comparison.
The comparison functions are based on those in
\cite{JianTu2023}. We use second differences for the boundary comparison, where
only the first derivatives extend continuously.

For $f$ in Proposition~\ref{extR:globalpair}, set $N_0=\sup_C f_{yy}\in(0,\infty)$. By the first comparison and the
strong maximum principle, equality holds for a normalized limit.
Comparing on the dual side, we obtain a constant horizontal vector
$\eta$ with $\eta\cdot e_y=1$ and $V_{\eta\eta}\le N_0^{-1}$.
At $P=Df(X)$, we have $1\le f_{yy}(X)V_{\eta\eta}(P)\le1$.
Equality implies $D^2f\,e_y=N_0\eta$. Integrating this identity and applying the same argument to the
other $y$ directions, we obtain the splitting.

\subsection{Notation}
For a convex function $u$, $Du$, $D^2u$ and $\partial u$
denote its gradient, Hessian and convex subdifferential. Its Alexandrov
measure is $M_u(E)=|\partial u(E)|$ for Borel sets $E$, where $|\cdot|$
is Lebesgue measure. We write $B_R$ for the Euclidean ball of radius
$R$ centered at zero, $\operatorname{cof}H$ for the cofactor matrix,
and $\operatorname{avg}_E$ for the Lebesgue average on $E$.

\section{Preliminaries}\label{sec:sections}\label{extR:chapter}
We establish estimates for centered sections, identify the gradient
image for the equation on a cone, and prove a local boundary estimate.

\subsection{Centered sections}
\label{extR:sections}

We define centered sections in the open effective domains of the minimal
convex extensions and use the mass balance identities to estimate their
Alexandrov measures.
When $D^2u>0$, write $(u^{ij})=(D^2u)^{-1}$, and similarly for the
other potentials. We sum repeated dummy coordinate indices only. Directional subscripts
denote differentiation along the specified constant vectors.

Fix $N\ge1$. Let $O,O'\subset\R^N$ be nonempty
open convex sets, and let
$\mu=g1_O\,dX$, $\nu=g'1_{O'}\,dP$ be locally finite measures,
with $g,g'>0$ a.e. on their respective supports.
Consider finite convex functions $u:O\to\R$, $v:O'\to\R$
whose restricted conjugates
\[
 \begin{gathered}
 U(X)=\sup_{p\in O'}\{X\cdot p-v(p)\},\\
 V(P)=\sup_{x\in O}\{x\cdot P-u(x)\}
 \end{gathered}
\]
satisfy $U|_O=u$, $V|_{O'}=v$, and the mass balance identities
\begin{equation}\label{extR:mass}
 \begin{gathered}
 \mu(E)=\nu(\partial U(E)),\\
 \nu(F)=\mu(\partial V(F)).
 \end{gathered}
\end{equation}
Here $\mathcal D_U=\operatorname{int}(\operatorname{dom}U)$,
$\operatorname{dom}U=\{U<\infty\}$, and similarly for $V$;
$E\subset\mathcal D_U$, $F\subset\mathcal D_V$ are Borel sets,
$\partial U(E)=\bigcup_{X\in E}\partial U(X)$, and measures are completed.
At boundary points we use the finite values of $U,V$. By definition,
$u(X)+V(P)\ge X\cdot P$ for $X\in O$, $P\in\R^N$, and
$U(X)+v(P)\ge X\cdot P$ for $X\in\R^N$, $P\in O'$.
These inequalities hold at finite boundary values by radial limits.
By separation from $\overline{O'}$, we have
$\partial U(\mathcal D_U)\subset\overline{O'}$.
Since $g'>0$ a.e. on $O'$ and $|\partial O'|=0$, \eqref{extR:mass}
implies that $M_U(E)=0$ if and only if $\mu(E)=0$.
Hence $\operatorname{supp}_{\mathcal D_U}M_U=\overline O\cap\mathcal D_U$
and $\operatorname{supp}_{\mathcal D_V}M_V=\overline{O'}\cap\mathcal D_V$.
Let $U^*(P)=\sup_{X\in\R^N}\{X\cdot P-U(X)\}$.
The ordinary convex conjugate of $U$ \cite{Rockafellar1970}
may differ from $V$ outside $O'$. For a finite transport plan
concentrated on $u(x)+v(p)=x\cdot p$, we have $U^*=v$ on $O'$
by restricted conjugacy. For a.e. $p\in O'$, we have $p\in\partial U(X)$
if and only if $X=Dv(p)$. Thus \eqref{extR:mass} follows from the
marginal identities \cite[Section~2 and Lemma~2.1]{FigalliJhaveri2023}.
For infinite measures we use change of variables. If the classical
interior gradients are inverse, the minimal extension is
\[
 U(X)=\sup_{x\in O}\{u(x)+Du(x)\cdot(X-x)\},
\]
with the analogous formula for $V$.

For $x_0\in\mathcal D_U$ and $H>0$, a complete centered section means
\[
 \begin{gathered}
 S_H^{c,U}(x_0)=\{X\in\mathcal D_U:
 U(X)<U(x_0)+p_H\cdot(X-x_0)+H\},\\
 \overline{S_H^{c,U}(x_0)}\Subset\mathcal D_U,\\
 \int_{S_H^{c,U}(x_0)}(X-x_0)\,dX=0.
 \end{gathered}
\]
We do not intersect the section with $O$. Its barycenter is taken
with respect to Lebesgue measure. The centers lie in
$\overline O\cap\mathcal D_U$. The slope $p_H$ need not be a
subgradient at $x_0$.
We write $S_H^c$ when the function is clear, and use the same
conventions for $V$.

The two density classes, invariant under invertible affine changes
and positive multiplication, are
\begin{enumerate}[label=(D\arabic*)]
\item $0<\lambda\le g\le\Lambda<\infty$ a.e. on $O$,
 with prescribed ratio $\frac{\Lambda}{\lambda}$;
\item $g=q^k$, where $k>0$ and $q:O\to(0,\infty)$ is finite and concave.
 Set $q=0$ outside $O$ when evaluating the density; concavity is required
 only on $O$.
\end{enumerate}
The two densities may belong to different classes with different constants.
Set $T_xE=x+\frac{E-x}{2}$. By concavity and interior approximation at
$x\in\partial O$, we have $q(T_xz)\ge \frac{q(z)}{2}$ for
$x\in\overline O$, $z\in O$.
Hence, for any Borel set $E$,
\begin{equation}\label{extR:doubling}
 \begin{gathered}
 \mu(E)\le C\mu(T_xE),\\
 C=\begin{cases}2^N\frac{\Lambda}{\lambda},&\text{(D1)},\\2^{N+k},&\text{(D2)}.
 \end{cases}
 \end{gathered}
\end{equation}
For these density classes, the mass identities imply a doubling estimate
for the Alexandrov measure on complete centered sections. We also estimate
the normalized potential near the section boundary.

\begin{lemma}\label{extR:modulus}
Suppose each density satisfies (D1) or (D2). Let
$S_H^{c,U}(x_0)$ be complete, $x_0\in\overline O$.
For $A\in GL(N)$ satisfying
$B_1\subset S=A^{-1}(S_H^{c,U}(x_0)-x_0)\subset B_D$, $D=N(N+1)$, set on
$A^{-1}(\mathcal D_U-x_0)$
\[
 w(X)=\frac{U(x_0+AX)-U(x_0)-p_H\cdot AX-H}{H}.
\]
Then
\[
 \begin{gathered}
 0\le-w(X)\le C\dist(X,\partial S)^\gamma,\\
 X\in S,
 \end{gathered}
\]
where $C$ depends only on dimension, exponents and density ratios,
and $\gamma=\frac{1}{N}$ when $g'$ satisfies (D1), or $\frac{1}{N+k}$ when
$g'$ satisfies (D2).
Moreover,
\[
 M_U\bigl(S_H^{c,U}(x_0)\bigr)
 \le C M_U\!\Bigl(x_0+\frac{1}{2}\bigl(S_H^{c,U}(x_0)-x_0\bigr)\Bigr).
\]
The same estimates hold for $V$.
\end{lemma}
\begin{proof}
We normalize the section by \cite[Lemma~2.4]{JS} and set
$\widetilde\mu(E)=\mu(x_0+AE)$ and
$\widetilde\nu(F)=\nu(p_H+HA^{-T}F)$.
These measures satisfy the same density classes and
$\widetilde\mu(E)=\widetilde\nu(\partial w(E))$.
Let $C_s$ be the doubling constant of $\widetilde\mu$, and denote
the transformed domain of $\widetilde\nu$ by $O'$. By convexity
along rays and the secant inequalities, we have $0\le-w\le M=1+D$
and $\partial w(\frac{1}{2}S)\subset B_L$, where $L=2M$.
For $K=\partial w(S)\subset\overline{O'}$ we have
\[
 \widetilde\nu(K)=\widetilde\mu(S)
 \le C_s\widetilde\mu(\frac{1}{2}S)\le C_s\widetilde\nu(B_L).
\]
Since $w(0)=-1$ and $w|_{\partial S}=0$, the minimum of
$w(x)-p\cdot x$ is attained in $S$ for $|p|<\frac{1}{D}$.
Hence $B_{\frac{1}{D}}\subset K$ and $0\in O'$.
For (D1), $\lambda_t|K|\le\widetilde\nu(K)$ and
$\widetilde\nu(B_L)\le\Lambda_t|B_L|$, since $|\partial O'|=0$.
For (D2), the density on $O'$ has the form $c_0q^{k_t}$.
Put $r=\frac{1}{2D}$ and $q_0=q(0)>0$.
By concavity on the segment from $z$ to $-\frac{rz}{|z|}\in O'$,
we have $q(z)\le\bigl(1+\frac{|z|}{r}\bigr)q_0$ for
$z\in O'\setminus\{0\}$.
At $z=0$ equality holds.
Moreover $\frac{1}{2}K\subset K$: for $p\in\partial w(x)$, $x\in S$,
\[
 (w-\frac{1}{2}p\cdot\mathrm{id})(Z)
 -(w-\frac{1}{2}p\cdot\mathrm{id})(x)\ge-\frac{w(x)}{2}>0
 \quad (Z\in\partial S),
\]
so the minimum is attained in the interior. By concavity and approximation at the endpoints,
$q\ge \frac{q_0}{2}$ on $\frac{1}{2}K\subset O'$. Therefore
\[
 2^{-N-k_t}c_0q_0^{k_t}|K|\le\widetilde\nu(\frac{1}{2}K)
 \le\widetilde\nu(K)\le C_s\Bigl(1+\frac{L}{r}\Bigr)^{k_t}c_0q_0^{k_t}|B_L|.
\]
Thus $|K|\le C$. In case (D1),
\[
 \lambda_t|B_{\frac{1}{D}}|
 \le\widetilde\nu(K)
 \le C_s\widetilde\nu\bigl(\partial w(\frac{1}{2}S)\bigr)
 \le C_s\Lambda_t|\partial w(\frac{1}{2}S)|.
\]
In case (D2), $B_{\frac{1}{2D}}\subset \frac{1}{2}K$ and $q\ge \frac{q_0}{2}$ there.
Since $\partial w(\frac{1}{2}S)\subset B_L$,
\[
 c\,c_0q_0^{k_t}
 \le\widetilde\nu(K)
 \le C_s\widetilde\nu\bigl(\partial w(\frac{1}{2}S)\bigr)
 \le C\,c_0q_0^{k_t}|\partial w(\frac{1}{2}S)|.
\]
Hence
\[
 c\le|\partial w(\frac{1}{2}S)|\le|\partial w(S)|\le C.
\]
It follows that $M_w(S)\le C M_w(\frac{1}{2}S)$.
The affine change of gradient multiplies both measures by
$\frac{H^N}{|\det A|}$, so the asserted inequality for $M_U$ follows.

The normalized extension is convex and lower semicontinuous.
Its negative set is $S$: an additional negative point on the boundary
of the effective domain would contradict completeness by convexity.
By \cite[Lemma~A.1]{CM} for $N\ge2$, and convex secants for $N=1$,
we have $(-w)^N\le C\dist(\cdot,\partial S)$.
For (D2), we multiply by $(-w)^{k_t}\le M^{k_t}$ to obtain the
stated exponent.
\end{proof}

The doubling estimate implies an engulfing property inside a complete
centered section.

\begin{lemma}[{\cite[Lemma~7.6]{CM}}]\label{extR:engulfing}
Let $u$ be finite and convex on an open convex set $\mathcal D$, and
set $X=\operatorname{supp}_{\mathcal D}M_u$.
Let $S=S_H^c(x)$ be complete, with $x\in X$ and $X\cap S$ convex.
Assume that, for all $y\in X\cap S$ and $0<b\le H$, the section
$S_b^c(y)$ is complete and
\[
 M_u\bigl(S_b^c(y)\bigr)
 \le C M_u\!\Bigl(y+\frac{1}{2}\bigl(S_b^c(y)-y\bigr)\Bigr).
\]
For $0\le t<\bar t\le1$ there exists $\sigma=\sigma(N,C,t,\bar t)>0$
such that
\[
 S_a^c(z)\subset x+\bar t(S-x)
\]
whenever $z\in X\cap\bigl(x+t(S-x)\bigr)$ and $0<a\le\sigma H$.
\end{lemma}
To apply \cite[Lemma~7.6]{CM}, extend $u$ by $+\infty$ outside
$\mathcal D$ and take its lower semicontinuous closure.
By convexity and completeness, no additional point below a section
plane lies on $\partial\mathcal D$. The complete sections and their
Alexandrov measures are therefore unchanged.

We apply the preceding lemma to both potentials to obtain strict
separation on their supports. We then use convex conjugacy to prove
gradient H\"older estimates, as in
\cite[Sections~4--6]{ForzaniMaldonado2004}.

\begin{lemma}\label{extR:holder}
Assume $U,V$ satisfy \eqref{extR:mass}, each density satisfies (D1) or
(D2), and complete centered sections exist at all heights and all
centers in the corresponding support. Then $U,V\in C^{1,\alpha}_{\loc}$ on their open effective
domains and are strictly convex on
$\overline O\cap\mathcal D_U$ and
$\overline{O'}\cap\mathcal D_V$.
If matching boundary points are normalized by $u(0)=v(0)=0$,
$0\in\partial U(0)\cap\partial V(0)$, and
$S_1^{c,U}(0),S_1^{c,V}(0)\subset B_K$, then $\mathcal D_U=\mathcal D_V=\R^N$ and, for all $R>0$,
\[
 \|U\|_{C^{1,\alpha}(B_R)}+\|V\|_{C^{1,\alpha}(B_R)}
 \le C(K,R,N,\text{density bounds}).
\]
\end{lemma}
\begin{proof}
By Lemma~\ref{extR:modulus}, $M_U,M_V$ satisfy the doubling inequality.
Their supports are convex by Section~\ref{extR:sections}, so
Lemma~\ref{extR:engulfing} applies to both potentials.

Subtract any affine support of $u|_{\overline O}$ at $x$, so that
$w\ge0$ on $\overline O$. At a boundary point, a supporting plane of the extension remains
a support of the restriction after adding an outward normal to its slope.
For $y\in\overline O\cap\partial S_H^c(x)$ choose $z=(1-t)x+ty$ with
$\frac{t}{1-t}>D$. Fix $t<t_1<t_2<1$. Since $\frac{t}{t_1}(y-x)+x$
lies in the open section, $z\in x+t_1(S_H^c(x)-x)$.
By Lemma~\ref{extR:engulfing} with $t_1,t_2$,
$S_{\sigma H}^c(z)\subset x+t_2(S_H^c(x)-x)$ and
$y\notin S_{\sigma H}^c(z)$.
If $x\in S_{\sigma H}^c(z)$, by balance we have
$z+\frac{z-x}{D}\in S_{\sigma H}^c(z)$.
Since $\frac{|z-x|}{D}>|y-z|$, convexity implies $y\in S_{\sigma H}^c(z)$,
a contradiction. The affine function $L$ defining this section satisfies
$L(x)\le0$ and $L(y)\le w(y)$. Since $z=(1-t)x+ty$,
\[
 \sigma H\le w(z)+\sigma H=L(z)=(1-t)L(x)+tL(y)\le t w(y).
\]
Hence
\begin{equation}\label{extR:gap}
 w(y)\ge cH.
\end{equation}
Applying \eqref{extR:gap} at the first intersection with
$\partial S_H^c(x)$ and using convexity, we obtain
$\{y\in\overline O:w(y)<cH\}\subset S_H^c(x)$.

Assume the normalization in the statement. For $|X|>K$, apply
\eqref{extR:gap} at the first intersection of the ray through $X$
with $\partial S_1^c(0)$. By convexity, $u(X)\ge \frac{c|X|}{K}$.
The same bound holds for $v$, and both functions are nonnegative.
Thus $U,V$ are bounded on $B_{2d}$, where $d>0$ depends only on
$K,N$, the density exponents and the density ratios.
If $1+\ell_1\cdot X$ defines $S_1^c(0)$, then
$1+\ell_1\cdot X-U(X)\le1+D$ on $B_{2d}$, so $|\ell_1|$ is bounded.
Since $U(0)=0$, by convexity we have $0\le U(X)\le \frac{C|X|}{d}$
on $B_d$. Hence, for $d'>0$ depending only on the same data,
$B_{d'}\subset S_1^c(0)$ for both potentials. Choose $\sigma$ small enough
that $S^c_{\sigma H}(0)\subset\frac{1}{2} S^c_H(0)$.
For each integer $j\ge1$, we obtain
\[
 \begin{gathered}
 S^c_{\sigma^{-(j-1)}}(0)\subset\frac{1}{2} S^c_{\sigma^{-j}}(0),
 \\
 2^jB_{d'}\subset2^jS_1^c(0)\subset S^c_{\sigma^{-j}}(0).
 \end{gathered}
\]
Thus both open effective domains are all of $\R^N$. For $H\ge1$,
\[
 d'|\ell_H|=\sup_{B_{d'}}\ell_H\cdot X\le DH+\sup_{B_{d'}}U\le CH.
\]
Choose the least $j$ such that $B_R\subset2^jB_{d'}$. Then, for some $a>1$,
we obtain
\[
 \sup_{B_R}U+\sup_{B_R}V\le C(1+R)^a .
\]
For $0\ne X\in\overline O$, take $P=\frac{tX}{|X|}$ in $u(X)\ge X\cdot P-V(P)$. Then
\[
 u(X)\ge\sup_{t\ge0}\{t|X|-C(1+t)^a\}.
\]
The same bound holds for $v$.

Put $G(X)=U(X)-\ell_1\cdot X$. By convexity along rays and the gap on the unit section, $G(X)\ge \frac{|X|}{K}-C$ on $\R^N$. If $Z\in T=S_1^c(x)$,
then $Z'=x-\frac{Z-x}{D}\in T$ by balance. Adding the inequalities
defining $T$ at $Z$ and $Z'$, we obtain
\[
 G(Z)+D\,G(Z')<(1+D)(G(x)+1).
\]
Thus $|Z|\le C$ for $x$ in each compact subset of the support. The same holds for $V$.
Iterating Lemma~\ref{extR:engulfing}, we obtain, for $h\le1$,
$\operatorname{diam}S_h^c(x)\le Ch^\beta$, where
$\beta=\frac{\log 2}{|\log\sigma|}$.
Combining this estimate with \eqref{extR:gap} and increasing $q=\frac{1}{\beta}$ above one if necessary, we obtain
\begin{equation}\label{extR:strict}
 u(x')-u(x)-p\cdot(x'-x)\ge c|x'-x|^q,
\end{equation}
for points of the support in each compact set and any slope $p$
supporting the restriction to the support at $x$, including boundary and
dual contacts. Choose $H$ comparable to
$|x'-x|^{\frac{1}{\beta}}$ with a sufficiently small constant so that $x'$ is outside
$S_H^c(x)$. Apply \eqref{extR:gap} at the first intersection of $[x,x']$ with the section boundary. For distances bounded below, fix one small $H$ and reduce $c$.
As the slope varies over a bounded set, all maximizers lie in one
compact subset of the closed support. Adding \eqref{extR:strict} for two maximizers, we obtain
\[
 c|x_1-x_2|^q\le(p_1-p_2)\cdot(x_1-x_2).
\]
Thus each maximizer is unique and depends H\"older continuously on
the slope, with exponent $\min\bigl(1,\frac{1}{q-1}\bigr)$. The two supporting
inequalities identify this maximizer with the gradient.
The argument also applies at interior centers. On $O'$, $U^*=v$:
we have $U^*(p)\le v(p)$, and for $z\in\partial v(p)$,
$U(z)=z\cdot p-v(p)$, so the reverse inequality holds.
For an interior ball $B\Subset O$, we have
$0<\mu(B)=\nu(\partial U(B))$, so we can choose $x\in B$ and
$p\in O'\cap\partial U(x)$. Then $U(x)+V(p)=x\cdot p$ and
$x\in\partial V(p)$. Translate $x,p$ to zero and subtract the
corresponding affine supports. Both sections of height one lie in
some $B_K$ by completeness. The preceding argument proves the
qualitative assertion.
\end{proof}

The preceding estimates require complete centered sections at all heights.
A linear lower bound after subtraction of an affine function, together
with divergence at the finite boundary, ensures their existence.

\begin{lemma}[{\cite[Theorem~A.7 and Lemma~A.8]{CM}}]\label{extR:centers}
Let $F$ be finite and convex on an open convex set $\mathcal D$.
Assume $F\to\infty$ at finite $\partial\mathcal D$ and
$F(x)-p_0\cdot x\ge c|x|-C$, with $c>0$.
For $x_0\in\mathcal D$ and $H>0$, there is a unique complete
centered section $S_H^c(x_0)\Subset\mathcal D$.
\end{lemma}
Extend $F$ by $+\infty$ outside $\mathcal D$. For $|p-p_0|<c$,
$F(x)-p\cdot x\ge(c-|p-p_0|)|x|-C$.
The extension is lower semicontinuous, and $F-p\cdot x$ attains its
minimum in $\mathcal D$. Hence $B_c(p_0)\subset\partial F(\R^N)$.
For $q\in\partial F(x_0)$, apply \cite[Theorem~A.7]{CM} to
\[
 \frac{F(x_0+z)-F(x_0)-q\cdot z}{H}.
\]
We obtain a bounded section centered at zero. Since $F\to\infty$
at finite boundary points, its translate is compactly contained in
$\mathcal D$. Uniqueness follows from \cite[Lemma~A.8]{CM}, including
uniqueness of the slope since
$\mathcal D\setminus\overline{S_H^c(x_0)}$ has nonempty interior.

By uniqueness, any invertible linear symmetry fixing the center and
preserving $F$ up to an affine function also preserves the section.

\subsection{The gradient image}
\label{extR:global}

On the cone below, the mass balance and the condition $f_r=0$ on the
boundary determine the gradient image. We can then apply
Lemma~\ref{extR:holder} to the minimal extensions.

\begin{proposition}\label{extR:globalpair}
Fix $k>0$. Let $h:\R^p\to[0,\infty)$ be finite, convex and
polyhedral, with $h(tz)=th(z)$ for $t\ge0$. Set
$C=\{r>h(z)\}\times\R_y$. Write $X=(r,z,y)$ and
$P=(s,\zeta,q)$, and set $N=p+2$. Let
$f\in C^\infty(C)\cap C^1_{\rm loc}(\overline C)$ satisfy
$f\ge0$, $D^2f>0$ and $f_r>0$ in $C$, with
\[
 \begin{gathered}
 f_r^k\det D^2f=K>0\quad\text{in }C,\\
 f_r=0\quad\text{on }\partial C.
 \end{gathered}
\]
Subtract $L_0(X)=f(0)+Df(0)\cdot X$ and still denote the result by $f$. Set
\[
 \begin{gathered}
 U(r,z,y)=f(\max(r,h(z)),z,y),\\
 V(P)=\sup_{X\in\overline C}(P\cdot X-f(X)).
 \end{gathered}
\]
Then $V<\infty$ on $\R^{p+2}$,
$Df:C\to\mathcal H=\{s>0\}$ is a diffeomorphism, and
$\det D^2V=\frac{s^k}{K}$ in $\mathcal H$.
The minimal extensions $U,V$ belong to $C^{1,\alpha}_{\loc}(\R^N)$
and are strictly convex on their respective closed supports.
\end{proposition}
\begin{proof}
Set $\mu=K\mathbf1_C\,dX$, $D=\{V<\infty\}$ and
$\nu=s^k\mathbf1_{D\cap\mathcal H}\,dP$.
By convexity and passage to the boundary, $L_0$ supports the
original potential on $\overline C$. Since $f_r(0)=0$, subtracting
$L_0$ preserves $f_r$, the equation and the boundary condition.
The normalized potential is nonnegative.

Suppose a sublevel of $w=f-\ell\ge0$ contains $Z_j$ with
$|Z_j|\to\infty$. After passing to a subsequence, we have
$\frac{Z_j}{|Z_j|}\to d\in\overline C$. By convexity on $[X,Z_j]$ and
passage to the limit, $w(X+td)\le w(X)$ for $t\ge0$.
Choose $B_{2\rho}(X_*)\Subset C$ and $T>4\rho$. Since $d$ is a
recession direction of $C$, we have $B_{2\rho}(X_*+jTd)\subset C$.
On these balls, $0\le w\le\sup_{B_{2\rho}(X_*)}w$.
By convexity, $Dw$, and hence $Df=Dw+D\ell$, are bounded by one
constant $L$ on the disjoint balls $B_j=B_\rho(X_*+jTd)$.
Since $Df$ is injective, by change of variables we obtain
\[
 K\sum_j|B_j|\le\int_{B_L\cap\mathcal H}s^k\,dP<\infty,
\]
a contradiction.
Thus $\{w\le w(X_*)+1\}$ is bounded for an interior point $X_*$.
To prove the linear bound, choose $R$ such that
$\{w\le w(X_*)+1\}\subset B_{\frac{R}{2}}(X_*)$. For $X\in C$ with $d=|X-X_*|\ge R$, put
$Z=X_*+\frac{R}{d}(X-X_*)$. By convexity,
\[
 w(X)\ge w(X_*)+\frac{d}{R}(w(Z)-w(X_*))\ge w(X_*)+\frac{d}{R}.
\]
Since $w\ge0$ for $d<R$, it follows that $w\ge c|X|-C'$ on $C$.
For $P\in\{V<\infty\}$, apply this to $f-P\cdot X+V(P)$:
we obtain $f(X)-P\cdot X\ge c_P|X|-C_P$. For $|P'-P|<\frac{c_P}{2}$,
$f(X)-P'\cdot X\ge\frac{c_P}{2}|X|-C_P$, whence $V(P')<\infty$.
Thus $D$ is open. By lower semicontinuity,
$V(P)\to+\infty$ at each finite point of $\partial D$.
At a contact point with $s>0$, the inward $r$ derivative of
$f-P\cdot X$ is $-s<0$. Therefore $Df(C)=D\cap\mathcal H$ and
$D^2V(Df(X))=(D^2f(X))^{-1}$.

Put $c(r,z,y)=(\max(r,h(z)),z,y)$. For $Y=c(X)$,
$a=(h(z)-r)_+$ and an affine support $L_Z$ at $Z\in C$,
\[
 L_Z(X)=L_Z(Y)-af_r(Z)\le f(Y).
\]
Since $af_r(Y)=0$, we let $Z\to Y$ from the interior and obtain
$U=f\circ c=\sup_{Z\in C}L_Z$. Hence $U$ is convex. Set
$G(X)=Df(c(X))$, using the continuous boundary trace. Then
$G(X)\in\partial U(X)$. Since $G$ is continuous on $\R^N$, we have
\[
 0\le U(X+Z)-U(X)-G(X)\cdot Z
 \le [G(X+Z)-G(X)]\cdot Z=o(|Z|).
\]
Thus $U\in C^1(\R^N)$ with $DU=G$.
By approximation with interior supports, $U$ is minimal and
\[
 V(P)=\sup_{Q\in D\cap\mathcal H}\{V(Q)+DV(Q)\cdot(P-Q)\}.
\]
For the extensions, we have
\[
 \begin{gathered}
 DU(\R^N\setminus C)\subset\{s=0\},\qquad \nu(\{s=0\})=0,\\
 \partial V(D\cap\{s\le0\})\subset\partial C,\qquad \mu(\partial C)=0.
 \end{gathered}
\]
The coercive bound for $f-P\cdot X$ is locally uniform in $P\in D$.
The set of maximizers
\[
 E(P)=\operatorname*{argmax}_{X\in\overline C}\{P\cdot X-f(X)\}
\]
is therefore nonempty, compact and convex, and $\partial V(P)=E(P)$ by the subdifferential formula for a supremum.
For $s\le0$, an interior maximizer would satisfy $s=f_r>0$.
Hence $E(P)\subset\partial C$, proving the second boundary inclusion.
Changing variables in the interior, we verify \eqref{extR:mass} for
$\nu$, whose density on $D\cap\mathcal H$ satisfies (D2).

For $P\in D\cap\mathcal H$, set $Y=c(X)$ and $a=(h(z)-r)_+$.
By coercivity,
\[
 U(X)-P\cdot X\ge c|Y|+sa-C\ge c'|X|-C.
\]
For $X\in C$, choose $B_\epsilon(X)\Subset C$. Since $f$ is bounded above on this ball,
\[
 V(P)-X\cdot P\ge \sup_{Y\in B_\epsilon(X)} P\cdot(Y-X)-C=\epsilon|P|-C.
\]
By Lemma~\ref{extR:centers}, complete centered sections exist for
$U,V$ at all centers in the support and all heights. Since $U,V\ge0$
and $U(0)=V(0)=0$, we have $0\in\partial U(0)\cap\partial V(0)$.
Choose $K_0$ such that $S_1^{c,U}(0)\cup S_1^{c,V}(0)\subset B_{K_0}$.
Applying Lemma~\ref{extR:holder} to the densities of types (D1) and
(D2), we obtain $D=\R^N$, the asserted regularity and strict convexity.
The other conclusions follow from the identities for $Df(C)$ and $D^2V(Df(X))$.
\end{proof}

Paired affine changes preserve the mass identities and doubling constants.
The resulting estimates are controlled by the radii of the two
normalized sections.

\begin{corollary}\label{extR:paired}
For a pair in Proposition~\ref{extR:globalpair},
let $x_0,p_0$ be matching boundary points. For $t>0$ and
$A\in GL(N)$ set
\[
 \begin{aligned}
 u_A(X)&=t^{-1}[u(x_0+AX)-u(x_0)-p_0\cdot AX],\\
 v_A(P)&=t^{-1}[v(p_0+tA^{-T}P)-v(p_0)-x_0\cdot tA^{-T}P].
 \end{aligned}
\]
Let $U_A,V_A$ be their minimal extensions. If the complete
centered sections satisfy $S_1^{c,U_A}(0)\subset B_K$ and
$S_1^{c,V_A}(0)\subset B_K$, then
\[
 \|U_A\|_{C^{1,\alpha}(B_R)}+\|V_A\|_{C^{1,\alpha}(B_R)}
 \le C(K,R,N,k).
\]
\end{corollary}
\begin{proof}
Using $u(x_0)+v(p_0)=x_0\cdot p_0$, we obtain the same formulas
for $U_A,V_A$ with $U,V$ in place of $u,v$.
The mass identities and doubling constants $2^N,2^{N+k}$ are
unchanged. We apply Lemma~\ref{extR:holder}.
\end{proof}

\begin{remark}
The estimate in Corollary~\ref{extR:paired} holds on all boundary
faces of the transformed domains. It also holds for a limiting pair
satisfying the same hypotheses, including the global gradient image
and boundary conditions.
\end{remark}

\subsection{Boundary estimates}
For local second derivative bounds, we use a direction along which
the boundary and the density are invariant. The following estimate
is independent of the lower bound for obliqueness.

\begin{lemma}\label{extR:directional}
Let $w$ be strictly convex and smooth up to the boundary of a
smooth convex domain $\Omega$, with
\[
 g'(Dw)\det D^2w=g.
\]
Let $\Omega'$ be the domain of $g'$. Assume $Dw(\Omega)\subset\Omega'$
and $Dw(X)\in\partial\Omega'$ at the points $X\in\partial\Omega$
under consideration. Let $K=\overline\Omega\cap\{w\le a\}$
be compact and separated from the edges of the coordinate neighborhood,
with $D^2w>0$ on $K$. Assume
\begin{enumerate}[label=(\roman*)]
\item $\partial\Omega$ is locally invariant under translations in the
$e_1$ direction near $K\cap\partial\Omega$, and $\partial_{e_1}g=0$
near $K$.
\item Near each corresponding point $P\in\partial\Omega'$, the domain
$\Omega'$ is a convex epigraph $P_N>\chi(P')$, with $\chi$
differentiable. Its inward
 conormal $b=(-D\chi,1)$ satisfies
$\chi(P')-D\chi(P')\cdot P'\le0$ and $b\cdot n>0$, where $n$ is the inward unit normal to $\partial\Omega$ at $X$,
with $P=Dw(X)$.
\item The interior densities are positive, and $D^2\log g'\le0$
and $P\cdot D\log g'(P)\le\kappa$.
\end{enumerate}
Then
\[
 (a-w)w_{11}\le C(N,\kappa,\|w_1\|_{L^\infty(K)})
 \quad\text{on }K,
\]
with $C$ independent of the positive obliqueness constant.
\end{lemma}
We assume the intercept condition in the original coordinates
$P=Dw(X)$, before translating the coordinate charts at $\partial\Omega'$.
\begin{proof}
For the interior calculation we follow \cite[Proposition~7.1]{CT}.
Let $\eta=a-w$, $\phi=\log g'$ and
$L=w^{ij}\partial_{ij}+D\phi(Dw)\cdot D$.
On $\{\eta>0\}$ set $\Psi=\log w_{11}+\log\eta+\frac{w_1^2}{2}$.
Let $\rho$ be the inward signed distance near $\partial\Omega\cap K$,
smoothly extended to $K$. At a boundary point with $\eta>0$, for $|t|$ small, we have
\[
 b\cdot\bigl(Dw(X+t e_1)-Dw(X)\bigr)\ge0.
\]
Differentiating at $t=0$ and using the intercept condition, we obtain
$b\cdot Dw_1=0$, $b\cdot Dw_{11}\ge0$, and
$b\cdot Dw=\chi(P')-D\chi(P')\cdot P'\le0$.
Hence $b\cdot D(\Psi+\varepsilon\rho)>0$ by $b\cdot n>0$.
Since $\Psi\to-\infty$ as $\eta\downarrow0$, the maximum is attained
at an interior point $X_\varepsilon$ whenever $\{\eta>0\}$ is nonempty.
We have $L\rho\ge-C$ on $K\cap\Omega$, where $C$ may depend on
$w,g',\Omega,\Omega'$ but is independent of $\varepsilon$.
The second order term is bounded on $K$. By strict obliqueness,
$P+tD\rho$ lies in a compact subset of $\Omega'$ near each boundary
contact, for one $t>0$ throughout that neighborhood. We bound
$\phi(P)$ above by a tangent plane and obtain
\[
 D\phi(P)\cdot D\rho(X)
 \ge\frac{\phi(P+tD\rho(X))-\phi(P)}{t}\ge-C.
\]
By compactness, finitely many such neighborhoods suffice.
Differentiating the logarithmic equation in the $e_1$ direction and
using $\partial_{e_1}g=0$, we obtain in the interior
\[
 \begin{aligned}
 Lw_1&=0,\\
 Lw_{11}&=w^{ia}w^{bj}w_{ab1}w_{ij1}
       -D^2\phi(Dw)[Dw_1,Dw_1],\\
 L\eta&=-N-Dw\cdot D\phi(Dw)\ge-(N+\kappa).
 \end{aligned}
\]
At $X_\varepsilon$, choose a constant frame with first vector $e_1$
such that $w_{ij}=\lambda_i\delta_{ij}$, $\lambda_1=w_{11}$ and
$\lambda_i=1$ for $i\ne1$.
By compactness and $D^2w>0$, we can bound the frame vectors, the first
derivatives in this frame and $\lambda_1^{\pm1}$ independently of
$\varepsilon$, with bounds depending on $w,\rho$. At $X_\varepsilon$, we have
\[
 \frac{w_{11i}}{\lambda_1}
 =\frac{w_i}{\eta}-w_1w_{1i}-\varepsilon\rho_i.
\]
The terms containing $w_{111}^2$ cancel. Combining the terms indexed
by $(1,i)$ and $(i,1)$, we obtain
\[
 \sum_{i,j}\frac{w_{ij1}^2}{\lambda_1\lambda_i\lambda_j}
 -\sum_i\frac{w_{11i}^2}{\lambda_i\lambda_1^2}
 \ge\sum_{i\ne1}\Bigl(\frac{w_{11i}}{\lambda_1}\Bigr)^2.
\]
The $-D^2\phi$ term is nonnegative. For $i\ne1$, we have $w_{1i}=0$.
By the first derivative identity,
\[
 \Bigl(\frac{w_{11i}}{\lambda_1}\Bigr)^2-\frac{w_i^2}{\eta^2}
 =\Bigl(\frac{w_i}{\eta}-\varepsilon\rho_i\Bigr)^2-\frac{w_i^2}{\eta^2}
 \ge-\frac{C\varepsilon}{\eta}.
\]
Substituting $L(\frac{w_1^2}{2})=\lambda_1$ and $L\rho\ge-C$, we obtain
\[
 0\ge L(\Psi+\varepsilon\rho)
 \ge \lambda_1-\frac{N+\kappa}{\eta}
       -\frac{w_1^2}{\lambda_1\eta^2}
       -C\varepsilon(1+\eta^{-1}).
\]
Set $Y_\varepsilon=\eta(X_\varepsilon)w_{11}(X_\varepsilon)$ and
$B=\|w_1\|_{L^\infty(K)}$.
Multiplying by $\lambda_1\eta^2$ and using
$\sup_K(\lambda_1+\eta)<\infty$, we obtain
\[
 Y_\varepsilon^2-(N+\kappa)Y_\varepsilon-B^2\le C\varepsilon .
\]
By the choice of $X_\varepsilon$, for $X\in K$ with $\eta(X)>0$,
\[
 \eta(X)w_{11}(X)
 \le e^{\frac{B^2}{2}+\varepsilon\operatorname{osc}_K\rho}\,Y_\varepsilon .
\]
Letting $\varepsilon\downarrow0$ proves the estimate.
\end{proof}

For families on an epigraph, we use bounded coordinate neighborhoods
that contain the sublevel sets and their gradient images at a uniform
positive distance from their boundaries.

\begin{proposition}\label{audit:general}
Fix $k,a>0$, let $\gamma:\R^{n-1}\to\R$ be finite and convex
with $\gamma(0)=0$, and set
$\Omega=\{r>\gamma(w)\}$, $\mathcal H=\{s>0\}$.
Consider a family of smooth strictly convex potentials $u,v$ for which
$Du:\Omega\to\mathcal H$ is a diffeomorphism and
\begin{equation}\label{extR:carrier-pair}
 \begin{gathered}
 V(P)=\sup_{X\in\overline\Omega}\{P\cdot X-u(X)\}<\infty
 \quad (P\in\R^n),\\
 u_r^k\det D^2u=K>0,
 \end{gathered}
\end{equation}
with $u_r=0$ on all of $\partial\Omega$ and $u,v\ge0$. Let the minimal extensions $U,V$ be locally $C^1$, with
$U(0)=V(0)=0$ and $DU(0)=DV(0)=0$.
Assume, uniformly in the family, $U,V\le M(R)$ on $B_R$ for $R\ge1$.
Choose $\epsilon,\tau$ independently of the member of the family, with
$B_\epsilon(\tau e_r)\subset\Omega$ and $u$ uniformly bounded
on this ball. Set $K_u=\overline\Omega\cap\{u\le4a\}$ and
$K_v=\overline{\mathcal H}\cap\{v\le4a\}$.
Let $\mathcal P,\mathcal P'$ be bounded open sets satisfying
\[
 \begin{gathered}
 \dist(K_u\cup DV(K_v),\R^n\setminus\mathcal P)\ge d_0>0,\\
 \dist(K_v\cup DU(K_u),\R^n\setminus\mathcal P')\ge d_0.
 \end{gathered}
\]
Fix a unit vector $\tau_0\in\R^{n-1}$. Suppose $\gamma$ is constant
on all segments parallel to $\tau_0$ in
$\{w:(\gamma(w),w)\in\mathcal P\}$.
Then
\[
 \begin{gathered}
 u_{\tau_0\tau_0}\le C\quad \hbox{on }\Omega\cap\{u<a\},\\
 v_{\eta\eta}\le C\quad \hbox{on }\mathcal H\cap\{v<a\},
 \quad \eta\perp e_s,\\
 |\eta|=1.
 \end{gathered}
\]
The constant depends only on $n,k,a$, the stated uniform bounds and
$d_0$, independently of $K$.
\end{proposition}
We prove this proposition after Proposition~\ref{loc:directional}.

\section{Section polarity and directional estimates}\label{sec:directional}
We first normalize the sections of $U$ and $V$ and obtain interior
estimates for $v_H$. We then compare second differences to bound
$(f_H)_{yy}$ on $A_H^{-1}C$ and $(v_H)_{ee}$ on $\mathcal H$ for
$e\perp e_s$.

\subsection{Sections and polar ellipsoids}
For an ellipsoid $E_H=A_HB_1$, write $E_H^\circ=A_H^{-T}B_1$.
We compare $S_H^{c,U}(0)$ with $E_H$ and $S_H^{c,V}(0)$ with
$HE_H^\circ$, using constants independent of $H$.
\begin{lemma}\label{extR:polarity}
For $U,V$ in Proposition~\ref{extR:globalpair} there exist
open ellipsoids $E_H=A_HB_1$ centered at zero, one for each $H>0$,
and constants $c,C>0$ independent of $H$, such that
$cE_H\subset S_H^{c,U}(0)\subset CE_H$ and the matrices $A_H$ have
the form
\[
 \begin{gathered}
 A_H=\begin{pmatrix}d_H&b_H^T&0\\0&M_H&0\\0&c_H^T&\rho_H\end{pmatrix},
 \\
 d_H,\rho_H>0,\qquad |b_H|\leq Cd_H,\\
 h(M_Hz)\leq Cd_H|z|\quad (z\in\R^p).
 \end{gathered}
\]
Moreover,
$cHE_H^\circ\subset S_H^{c,V}(0)\subset CHE_H^\circ$.
There are
$0<\beta\leq\frac{1}{2}$, $\alpha=1-\beta<1$ such that
\[
 c\Bigl(\frac{R}{H}\Bigr)^\beta E_H\subset E_R\subset C\Bigl(\frac{R}{H}\Bigr)^\alpha E_H\quad (R\geq H).
\]
The normalized extensions are
$U_H(X)=H^{-1}U(A_HX)$ and $V_H(P)=H^{-1}V(HA_H^{-T}P)$.
Set
$f_H=U_H|_{A_H^{-1}C}$ and $v_H=V_H|_{\mathcal H}$.
The extensions $U_H,V_H$ have polynomial upper bounds on $\R^N$
independent of $H$. Both $f_H$ and $v_H$ have a superlinear lower
bound independent of $H$. The domains $A_H^{-1}C$ contain an interior
ball independent of $H$, and their graph functions have Lipschitz
constants bounded independently of $H$.
\end{lemma}
\begin{proof}
Choose a centered comparison ellipsoid satisfying
$E\subset S_H^c(0)\subset C_NE$. Taking the Cholesky factorization in the order $(z,r,y)$, we write
\[
 \begin{gathered}
 E=A B_1,\\
 A=\begin{pmatrix}d&b^T&0\\0&M&0\\a&c^T&\rho\end{pmatrix},
 \\
 d,\rho>0.
 \end{gathered}
\]
Let $p_r$ be the $r$ component of the slope of the section plane.
We have $p_r>0$. Otherwise, the negative $r$ axis lies in the section,
since $U=0$ on this axis.
For $c_0(r,z,y)=(\max\{r,h(z)\},z,y)$, the difference between the
section plane and $U$ increases by $p_r(h(z)-r)_+$. Thus
$c_0(\overline E)\subset C_N\overline E$.
Applying this inclusion at $A(-1,0,0)=(-d,0,-a)$, we obtain
$A^{-1}c_0(A(-1,0,0))=(0,0,-\frac{a}{\rho})$ and $|a|\le C_N\rho$.
At $A(0,z,0)$, $|z|=1$, the first component is
$\frac{(h(Mz)-b\cdot z)_+}{d}\le C_N$. Applying this also to $-z$ and using
$h\ge0$, we have
$|b\cdot z|\le C_Nd$ and $h(Mz)\le2C_Nd$.
Let $S(t,z,y)=(t,z,y-\frac{at}{\rho})$. The first column of $AS$
is $(d,0,0)^T$. Both $\|S\|$ and $\|S^{-1}\|$ are at most
$L=1+C_N$ because $\frac{|a|}{\rho}\le C_N$. Thus
$L^{-1}AB_1\subset ASB_1\subset L AB_1$.
Set $A_H=AS$ and $E_H=ASB_1$. The comparison constants are independent
of the height. By subadditivity, $\frac{h(Mz)-b\cdot z}{d}$ is uniformly Lipschitz.
Its epigraph contains a ball centered on the positive $r$ axis, with
center and radius independent of $H$.

Suppose $B_d\subset A_H^{-1}S_H^{c,U}(0)$ and the normalized section
plane is $1+m\cdot X$. For $|e|=1$, we have
$0\le U_H(\pm\frac{d}{2}e)<1\pm\frac{d}{2}m\cdot e$.
Hence $|m|\le\frac{2}{d}$.
By the Lipschitz bound for the graph, we can choose
$\tau,\epsilon>0$ sufficiently small so that
$B_\epsilon(\tau e_r)\subset B_d\cap A_H^{-1}C$. On this ball $f_H=U_H\le1+d|m|$. By the definition of $V_H$,
\[
 V_H(P)\ge \sup_{X\in B_\epsilon(\tau e_r)}P\cdot X-C=\tau s+\epsilon|P|-C.
\]
If $P$ belongs to the dual section $T$ of height one, then
$-\frac{P}{D}\in T$ by balance. Adding the inequalities defining $T$ at
these two points, we obtain
$V_H(P)+D V_H(-\frac{P}{D})\le1+D$.
Both terms are nonnegative. At least one of the two points has
$s\ge0$ and norm at least $\frac{|P|}{D}$. Hence $T\subset B_C$.
Both sections of height one have bounded outer radii. As in the
proof of Lemma~\ref{extR:holder}, we obtain
$B_c\subset T\subset B_C$.

For $0<a\le H$, normalize $S_H^c$ and let $w$ be the difference
between the potential and its section plane. By convexity along rays,
$w(0)=-H$ and $w(X)\ge cH|X|-CH$.
Apply the estimate at opposite points from the proof of
Lemma~\ref{extR:holder} to the transformed $S_a^c$, whose center
value is $-H$. We obtain
$S_a^c\subset CS_H^c\quad(0<a\le H)$.

By Lemma~\ref{extR:engulfing}, choose
$\theta\in(0,1)$ sufficiently small relative to the two ellipsoid
comparison constants that
$2E_{\theta H}\subset E_H$.
For $j\ge0$ with $\theta^{-j}\le \frac{R}{H}<\theta^{-j-1}$, iteration and
$S_a^c\subset CS_H^c$ imply
$2^jE_H\subset E_{\theta^{-j}H}\subset CE_R$.
Set $\beta_U=\frac{\log 2}{\log\bigl(\frac{1}{\theta}\bigr)}$. Since
$2^j>\frac{1}{2}\bigl(\frac{R}{H}\bigr)^{\beta_U}$, this proves the lower inclusion.
Apply the same argument to $V$, using $E'_h=hE_h^\circ$. For some
$\beta_V>0$,
$c\bigl(\frac{R}{H}\bigr)^{\beta_V}HE_H^\circ\subset RE_R^\circ$.
Taking polar sets, we obtain
$E_R\subset c^{-1}\bigl(\frac{R}{H}\bigr)^{1-\beta_V}E_H$.
Choosing $\beta=\min(\beta_U,\beta_V,\frac{1}{2})$ proves both bounds.

By the iteration and Fenchel estimate in the proof of
Lemma~\ref{extR:holder}, the polynomial upper bounds for $U_H,V_H$
and the superlinear lower bounds for $f_H,v_H$ are independent of $H$.
\end{proof}

\subsection{Interior estimates and comparison}
The domains $C$ and $A_t^{-1}C$ are invariant under translation in
the direction $e_y$. We first obtain interior estimates for $v_t$
on $\mathcal H$ and use the inverse gradient map to estimate $f_t$.

\begin{lemma}\label{audit:interior}
Let $f_t$ and $v_t$, $t>0$, be as in Lemma~\ref{extR:polarity}. For each $K\Subset\mathcal H$ there are
$c_K,C_K,d_K>0$, independent of $t$, such that, for all $P\in K$,
\[
 \begin{gathered}
 c_KI\le D^2v_t(P)\le C_KI,\\
 B_{d_K}(Dv_t(P))\subset A_t^{-1}C,\\
 D^2f_t\le C_KI\quad \hbox{on }B_{d_K}(Dv_t(P)).
 \end{gathered}
\]
For each integer $m\ge1$, the derivatives $D^mv_t$ are uniformly
bounded on $K$.
\end{lemma}
\begin{proof}
Write $N=p+2$ and omit the subscript $t$ from $f_t,v_t$ in this proof.
Since $A_te_r=d_te_r$, by the chain rule we have
\[
 \begin{gathered}
 f_r^k\det D^2f=\kappa,\\
 \det D^2v=\frac{s^k}{\kappa},\\
 \kappa=\frac{K_0d_t^k(\det A_t)^2}{t^{N+k}},
 \end{gathered}
\]
where $K_0$ is the original equation constant. Choose the ball
$B\subset A_t^{-1}C$ from Lemma~\ref{extR:polarity}, independently
of $t$, and choose a ball $B'\Subset\mathcal H$ independent of $t$. By the bounds on the convex extensions and
convexity, we have
$Df(B)\cup Dv(B')\subset B_C$.
Changing variables on $B$ and $B'$, we obtain
\[
 \begin{gathered}
 \kappa |B|=\int_{Df(B)}s^k\,dP\le C,\\
 \kappa^{-1}\int_{B'}s^k\,dP=|Dv(B')|\le C.
 \end{gathered}
\]
Hence $c\le\kappa\le C$.

Fix $K\Subset\mathcal H$ and $d>0$ such that
$\{Q:\dist(Q,K)\le3d\}\subset\mathcal H$.
For $P\in K$, set
$w(Q)=v(Q)-v(P)-Dv(P)\cdot(Q-P)$.
By Lemma~\ref{extR:holder}, $w\ge cd^q$ on $\partial B_d(P)$.
Set $h=\frac{cd^q}{2}$, $S=\{w<h\}$ and $z=w-h$. Then
$\overline S\Subset B_d(P)$, and $z$ is smooth and strictly convex, with
$z=0\quad\text{on }\partial S$, $\|z\|_\infty=h$ and $\|Dz\|_\infty\le C$.
On the larger compact set, $\frac{s^k}{\kappa}$ has uniform positive upper
and lower bounds and a uniform $C^2$ bound. Since $z(P)=-h$, the
Pogorelov estimate \cite[Theorem~2.12, pp.~539--540]{DPF} implies
$h^2w_{ee}(P)\le C\quad(|e|=1)$.
Thus $D^2v\le CI$ on $K$, since $h$ is independent of $t$. For the lower bound,
\[
 \lambda_{\min}(D^2v)\ge\frac{\det D^2v}{\lambda_{\max}(D^2v)^{N-1}}\ge c\quad\text{on }K.
\]
For higher regularity, choose $m,M>0$ so that $mI\le D^2v\le MI$
on a larger compact subset of $\mathcal H$. For symmetric $A$, set
\[
 E(A)=\inf_{\frac{m}{2}I\le B\le2MI}
 \Bigl\{\log\det B+\tr\bigl(B^{-1}(A-B)\bigr)\Bigr\}.
\]
By concavity of $\log\det$,
$E(A)=\log\det A\quad\bigl(\frac{m}{2}I\le A\le2MI\bigr)$.
As an infimum of affine functions, $E$ is concave.
For symmetric $A$ and $Q\ge0$,
$\frac{\tr Q}{2M} \le E(A+Q)-E(A) \le\frac{2\tr Q}{m}$.
Thus $E$ is uniformly elliptic on the space of symmetric matrices and
$E(D^2v)=g=k\log s-\log\kappa$.
Since $g$ has uniformly bounded derivatives of each order on this compact set,
we apply the interior Evans--Krylov estimate
\cite[Theorem~3.4 and Lemma~4.1]{TosattiWangWeinkoveYang2015} to obtain
$\|v\|_{C^{2,\alpha}(K)}\le C_K$
for some $\alpha\in(0,1)$, uniformly in $t$.
The equation for $v_a$ is
$v^{ij}(v_a)_{ij}=g_a$.
Applying the interior Schauder estimate
\cite[Theorem~6.2]{GilbargTrudinger2001} to $v_a$ and differentiating
the equation, we obtain uniform bounds for $D^mv$ on smaller compact
sets by induction on $m$.

For $P\in K$ take $\rho>0$ with $D^2v\ge cI$ on $B_{2\rho}(P)$,
uniformly in $P$. If $|Y-Dv(P)|<\frac{c\rho}{2}$, then on
$\partial B_\rho(P)$,
\[
 v(Q)-v(P)-Y\cdot(Q-P)\ge \frac{c\rho^2}{2}-|Y-Dv(P)|\rho>0.
\]
The minimum of $v(Q)-Y\cdot Q$ on $\overline B_\rho(P)$ is attained
in the interior. Hence
$B_{\frac{c\rho}{2}}(Dv(P))\subset Dv(B_\rho(P))\subset A_t^{-1}C$.
Using $D^2f(Dv(Q))=(D^2v(Q))^{-1}$, we obtain the asserted bound
for $D^2f$.
\end{proof}

We now compare second differences of $f$ and $v$ on sublevel sets.
Use $f$ from Proposition~\ref{extR:globalpair} and
$v=V|_{\mathcal H}$, or use $f_t,v_t$ from
Lemma~\ref{extR:polarity}, omitting the subscript $t$. Put $\ell=r$ in the original coordinates
and $\ell=r+\frac{b_t}{d_t}\cdot z$ in the normalized coordinates.
Write $\omega=D\ell$ and set
\[
 (g,J,a,e)=
 \begin{cases}
 (f,\ell f_r,\frac{N+k}{k+2},e_y),\\
 (v,sD_\omega v,\frac{N}{k+2},e),&e\perp e_s.
 \end{cases}
\]
For $h>0$, $0<\delta<\frac{1}{a}$ and $\tau\ne0$, define
$F_h=h+g-aJ$, $O_h=\{J<\delta(h+g)\}$,
$\Sigma_h=\{J=\delta(h+g)\}$ and
\[
 w=\frac{g(\,\cdot+\tau e)+g(\,\cdot-\tau e)-2g}{\tau^2}.
\]
\begin{lemma}\label{loc:front-comparison}
Assume $\frac{w}{h+g}\to0$ as $g\to\infty$ in $\overline O_h$.
Let $m\in\R$ and $\lambda\ge0$. If $w\le m+\lambda F_h$ on
$\Sigma_h$, then
$w\le m+\lambda F_h\quad\text{in }O_h$.
\end{lemma}
\begin{proof}
We have $\ell\ge0$ on $\overline C$ or $\overline{A_t^{-1}C}$,
respectively.
With $L_f=f^{ij}\partial_{ij}+\frac{k}{f_r}\partial_r$ and
$L_v=v^{ij}\partial_{ij}$, we differentiate the equations and obtain
\[
 \begin{gathered}
 L_f f=N+k,\\
 L_f(\ell f_r)=k+2,\\
 L_vv=N,\\
 L_v(sD_\omega v)=k+2.
 \end{gathered}
\]
Thus $L_gF_h=0$. By concavity of
$\log\det D^2f+k\log f_r$ and $\log\det D^2v$, respectively,
we have $L_gw\ge0$. The translations in directions $e\perp e_s$ preserve $s^k$.
Also $F_h\ge(1-a\delta)(h+g)>0$ on $\overline O_h$.

At $(r_0,z,y)$ on the boundary of the domain of $f$, we have $w_r=0$.
For $d>0$, by convexity and $a\ge1$,
\[
 F_h(r_0+d,z,y)-F_h(r_0,z,y)
 \le[-a\ell(r_0,z,y)-(a-1)d]f_r(r_0+d,z,y)\le0.
\]
For the boundary calculation for $v$, write $Z=(z,y)$ and the domain
of $f$ as $r>\gamma(Z)$, where $\gamma$ is finite and convex with
$\gamma(0)=0$. At $P\in\partial\mathcal H$, let
$X=Dv(P)=(\gamma(Z_0),Z_0)$ be its unique contact point.
For $b\in\partial\gamma(Z_0)$, set
$\nu=(1,-b)$ and
\[
 c=\nu\cdot X=\gamma(Z_0)-b\cdot Z_0\le\gamma(0)=0.
\]
The supporting plane inequality implies $\nu\cdot Y\ge c$ on
$\{r\ge\gamma(Z)\}$. In particular,
$\nu\cdot Dv(P\pm\tau e)\ge c$.
Using $\nu_s=1$ and the $C^1$ traces, we obtain
\[
 \begin{gathered}
 D_\nu w(P)\ge0,\\
 D_\nu^+F_h(P)=c-aD_\omega v(P)=c-a\ell(X)\le0.
 \end{gathered}
\]

Fix $\epsilon>0$. For all sufficiently large $R$, the assumptions on $\Sigma_h$
and at infinity imply that
$w-m-(\lambda+\epsilon)F_h$ is strictly negative on $\Sigma_h$ and
on $\{g=R\}$ of $O_h\cap\{g<R\}$. The closure of $O_h\cap\{g<R\}$ is compact. Therefore, for all sufficiently small $\mu>0$, adding
$\mu r$ when $g=f$, or $\mu(s+s^2)$ when $g=v$, preserves those
strict inequalities. This function cannot have a positive
interior maximum: its linearized operator is at least $\frac{\mu k}{f_r}>0$
or $2\mu v^{ss}>0$, respectively.

On the boundary of the domain of $g$, we have $J=0<\delta(h+g)$,
so a short inward segment lies in the truncation. The trace inequalities
imply that the inward lower derivative is at least $\mu$. Thus a
positive maximum is impossible there.
Let $\mu\downarrow0$ at each fixed $\epsilon,R,\tau$, then let $R\to\infty$, and finally let $\epsilon\downarrow0$.
\end{proof}

We combine Lemmas~\ref{audit:interior} and~\ref{loc:front-comparison}
to bound $(f_t)_{yy}$ on $A_t^{-1}C$ and $(v_t)_{ee}$ on
$\mathcal H$ for $e\perp e_s$.

\begin{proposition}\label{loc:directional}
Let $f_t$ and $v_t$ be as in Lemma~\ref{audit:interior}. Then $(f_t)_{yy}\le C(1+f_t)$ in $A_t^{-1}C$ and
$(v_t)_{ee}\le C(1+v_t)$ in $\mathcal H$, for all horizontal unit
vectors $e$ in the second inequality, with $C$ independent of $t$.
\end{proposition}
\begin{proof}
For $t>0$, omit the subscript from $f_t,v_t$ and set
$\ell_t(r,z,y)=r+B_t\cdot z$, $B_t=\frac{b_t}{d_t}$ and $\omega_t=D\ell_t=(1,B_t,0)$.
Since $d_t\ell_t$ is the original $r$ coordinate, we have
$\ell_t\ge0$ on $\overline{A_t^{-1}C}$. The vectors $\omega_t$
are uniformly bounded. Set
$J=\ell_t f_r$ and $I=sD_{\omega_t}v$.
They are nonnegative and agree when $P=Df(X)$.
For a second height $Rt$, $R\ge1$, define
$T=A_t^{-1}A_{Rt}$ and $T'=RT^{-T}$. From the two ellipsoid inclusions we have
$\|T^{-1}\|\le CR^{-\beta}$ and $\|(T')^{-1}\|=\frac{\|T\|}{R}\le CR^{-\beta}$.
We have $Te_r=\frac{d_{Rt}}{d_t}e_r$ and
$\ell_t(TY)=\frac{d_{Rt}}{d_t}\ell_{Rt}(Y)$. With the functions at height
$Rt$,
$\tilde f(Y)=\frac{f(TY)}{R}$ and $\tilde v(Q)=\frac{v(T'Q)}{R}$,
we obtain
\[
 \begin{gathered}
 R^{-1}J(TY)=\ell_{Rt}(Y)\tilde f_r(Y),\\
 R^{-1}I(T'Q)=Q_sD_{\omega_{Rt}}\tilde v(Q).
 \end{gathered}
\]
The second identity follows at $Y=D\tilde v(Q)$.

Set $a_0=\frac{N+k}{k+2}\ge1$ and $c_0=\frac{N}{k+2}$. Fix $\delta>0$
with $\delta\max(a_0,c_0)<1$. Take $R=1+f(X)$ for
$X\in A_t^{-1}C$, and $R=1+v(P)$ for $P\in\mathcal H$. The normalized value is $\frac{R-1}{R}\le1$.
By the superlinear lower bounds, the normalized points are bounded.
Using convexity and the local upper bounds, we also bound their gradients.
If $J(X)\ge\delta R$, then $Y=T^{-1}X$ and
$Q=D\tilde f(Y)$ satisfy $Q_s\ell_{Rt}(Y)\ge\delta$.
It follows that $Q_s\ge c>0$. By Lemma~\ref{audit:interior}, $D^2\tilde f$ is uniformly bounded above on a ball about $Y$ with a uniform radius. Likewise $I(P)\ge\delta R$ implies
$Q_sD_{\omega_{Rt}}\tilde v(Q)\ge\delta$ for $Q=(T')^{-1}P$,
hence $Q_s\ge c$ and $D^2\tilde v$ is bounded on a ball about $Q$ with a uniform radius.

For $g=f$ or $g=v$ and a unit vector $e$, set
\[
 \Delta_{\tau,e}^2g(X)
 =\frac{g(X+\tau e)+g(X-\tau e)-2g(X)}{\tau^2}.
\]
Take $e=e_y$ for $f$, and an arbitrary $e\perp e_s$ for $v$.
Choose $\tau_0>0$ sufficiently small, independently of $t,R,e$.
For $0<|\tau|\le\tau_0$, the integration segments below lie in these
balls. For $f$,
\[
 \Delta_{\tau,e}^2f(X)
 =R\int_{-1}^1(1-|a|)
 D^2\tilde f(Y+a\tau T^{-1}e)[T^{-1}e,T^{-1}e]\,da ;
\]
For $v$, use $T'$ in place of $T$. Thus each second difference
is bounded by $C R^{1-2\beta}$ where $J\ge\delta R$ or
$I\ge\delta R$, respectively. On the whole domain, by the Lipschitz bound on the normalized ball,
\[
 \begin{gathered}
 0\le\Delta_{\tau,e}^2f(X)\le\frac{2CR^{1-\beta}}{|\tau|},\\
 0\le\Delta_{\tau,e}^2v(P)\le\frac{2CR^{1-\beta}}{|\tau|}
 \end{gathered}
\]
For fixed $\tau\ne0$ these bounds are $o(R)$. The differences extend continuously to the boundary since the translations preserve the domain.

On the two regions $O_f=\{J<\delta(1+f)\}$ and
$O_v=\{I<\delta(1+v)\}$ put
$F_f=1+f-a_0J$ and $F_v=1+v-c_0I$.
On these sets, $\frac{F_f}{R}\ge1-a_0\delta>0$ and
$\frac{F_v}{R}\ge1-c_0\delta>0$, respectively. On their defining level sets, the $CR^{1-2\beta}$ estimates imply
$\Delta_{\tau,e}^2g\le C_0F_g$, with $C_0$ independent of
$t,\tau$ and the unit direction $e\perp e_s$. Using the $o(R)$ bounds, we apply Lemma~\ref{loc:front-comparison} with $h=1$, $m=0$ and $\lambda=C_0$ to extend both inequalities to $O_f,O_v$, respectively.
Since $0<F_g\le R$ on these regions and $CR^{1-2\beta}\le CR$
on their complements, we obtain
$\Delta_{\tau,e}^2g\le CR$
on $A_t^{-1}C$ when $g=f$ and on $\mathcal H$ when $g=v$.
At an interior point, let $\tau\to0$ to obtain the Hessian bounds.
\end{proof}

\begin{proof}[Proof of Proposition~\ref{audit:general}]
Write $X=(r,z)$, $P=(s,p)$ and $e=\tau_0$, with
$z,p\in\R^{n-1}$. For translation in the direction $d$, set
\[
 Q_h[b](Z)=\frac{b(Z+hd)+b(Z-hd)-2b(Z)}{h^2}.
\]

For $a_*\in\partial\gamma(0)$ and $|z|<\epsilon$, we have
$(\tau,z)\in\Omega$ and $a_*\cdot z\le\gamma(z)<\tau$. Thus $|a_*|\le\frac{\tau}{\epsilon}$. Since $\gamma$ is locally constant
along $e$, we also have $a_*\cdot e=0$.
The bounded shear $r_{\rm old}=r+a_*\cdot z$ therefore preserves
$e$ and replaces $\gamma$ by $\gamma-a_*\cdot z\ge0$.
Its dual map preserves $s$ and the subspace $s=0$, so these bounds
and all tangential directional estimates remain uniform.

As in Proposition~\ref{extR:globalpair}, $U(r,z)=u(Y)$ for
$Y=(\max\{r,\gamma(z)\},z)$. By conjugacy,
$u(X)\ge L|X|-M(L)$ and $v(P)\ge L|P|-M(L)$ on their respective
domains. Consequently
\[
 U(X)-r=u(Y)-Y_r+(Y_r-r)\ge |Y|+(Y_r-r)-M(2)\ge |X|-M(2).
\]

On $B_{\epsilon'}(\tau e_r)$ we have $u\le C$. Hence
\[
 V(P)-\tau s\ge \sup_{X\in B_{\epsilon'}(\tau e_r)}P\cdot(X-\tau e_r)-C=\epsilon'|P|-C.
\]
The centered sections exist by Lemma~\ref{extR:centers}. Applying
the affine lower bounds at opposite points as in
Lemma~\ref{extR:holder}, we bound both unit sections. Changing
variables as in Proposition~\ref{extR:globalpair}, we obtain
\eqref{extR:mass}. Hence \eqref{extR:strict} holds uniformly on
bounded sets. Also
$0\le\gamma(z)\le U(0,z)+M(1)$, so $\gamma$ is uniformly Lipschitz
on compact subsets of $\R^{n-1}$, uniformly in the family.
\smallskip\noindent\emph{Estimates for $v$ on $\mathcal H$.}
For $H>0$, choose a comparison ellipsoid with
\[
 \begin{gathered}
 A_H=\begin{pmatrix}d_H&b_H^T\\0&M_H\end{pmatrix},\\
 d_H>0.
 \end{gathered}
\]
Using $r\mapsto\max\{r,\gamma(z)\}$ as in
Lemma~\ref{extR:polarity}, we obtain $|b_H|\le Cd_H$ and
$\gamma(M_Hz)\le Cd_H\quad(|z|=1)$.
Thus the normalized graph
$\gamma_H(z)=\frac{\gamma(M_Hz)-b_H\cdot z}{d_H}$
is uniformly bounded above on $B_1$, bounded below by $-C|z|$,
and uniformly Lipschitz on $B_{\frac{1}{2}}$.
Since $\gamma_H(0)=0$, $A_H^{-1}\Omega$ contains a ball centered
on the positive $r$ axis inside the unit section. Its center and radius
can be chosen independently of $H$ and of $(u,v)$.
Applying the proof of Lemma~\ref{extR:polarity}, we obtain
$\|A_1\|+\|A_1^{-1}\|\le C$.

Apply the dual argument in Proposition~\ref{loc:directional} using
these bounds and Lemma~\ref{audit:interior}. At
$P\in\partial\mathcal H$, use a supporting plane at
$Dv_H(P)\in\partial(A_H^{-1}\Omega)$. The sign condition follows from
$r+\frac{b_H}{d_H}\cdot z=\frac{r_{\rm old}}{d_H}\ge0$.
Thus $v_{\eta\eta}$ is uniformly bounded above on
$\mathcal H\cap\{v<a\}$ for each $a>0$. Since $Du$ and $v(Du)$
are uniformly bounded on bounded subsets of $\Omega$, we have
\[
 (u^{\alpha\beta})_{1\le\alpha,\beta\le n-1}=(D^2_{pp}v)(Du)\le C I . \tag{A}
\]
Integrating as in the proof of Lemma~\ref{audit:interior}, we bound
$K$ above and below. By (A),
$u_{rr}=(\det D^2u)\det(u^{\alpha\beta})\le Cs^{-k}$, where $s=u_r$.
Set $d(X)=r-\gamma(z)$. Since $s=0$ on the graph,
\[
 s^{k+1}(r,z)=\lim_{\varepsilon\downarrow0}\int_{\gamma(z)+\varepsilon}^{r}\partial_r(s^{k+1})(t,z)\,dt\le C d(X).
\]
On each bounded subset of $\overline\Omega$, the uniform local
Lipschitz bound for $\gamma$ implies $\dist(X,\partial\Omega)\ge c d(X)$.
In the following estimates, the constants may depend on this subset.

At $X_0$, set $w(X)=u(X)-u(X_0)-Du(X_0)\cdot(X-X_0)$.
Let $s_0=u_r(X_0)>0$. Choose $c'>0$, uniformly in the family, so that
$\overline B_{2R}(X_0)\subset\Omega$ for $R=c'\min\{1,d(X_0)\}$.
By uniform support separation, $w\ge a_0R^q$ on $\partial B_R(X_0)$
with uniform constants $a_0>0$ and $q\ge2$. By convexity, for
$h=\frac{a_0R^q}{2}$ we have $\overline{\{w<h\}}\subset B_R(X_0)$.

In the translated gradient coordinates $Q=(\sigma,\pi)$, let
$\phi(Q)=k\log(\sigma+s_0)$. We have $D^2\phi\le0$ and
$Q\cdot D\phi=\frac{k\sigma}{\sigma+s_0}\le k$.
Since $|Dw|$ is bounded and $h-w=0$ on the section boundary, we apply
the interior calculation in Lemma~\ref{extR:directional} with cutoff
$h-w$ to obtain $h u_{\xi\xi}(X_0)\le C\quad(|\xi|=1)$.
Hence $D^2u(X)\le C(1+d(X)^{-q})I$.

Fix $X_0=(r_*,z_*)\in\Omega\cap\{u<a\}$. We have $|X_0|\le C$.
Since $u_r>0$, $u(\gamma(z_*),z_*)\le u(X_0)<a$.
Choose $\ell,\rho>0$ uniformly so that the graph above
$B_{2\rho}(z_*)$ lies in $\mathcal P$ and
$r_T=\gamma(z_*)+\ell\ge\gamma(z)+\frac{\ell}{2}\quad(z\in B_{2\rho}(z_*))$.
It suffices to consider $r_*<r_T$. The graph is constant along $e$
in this ball. On $r=r_T$, \eqref{extR:strict} implies
$c d^q\le u(r,z)-u(\gamma(z),z)\le u_r d$. Hence $u_r\ge c_\ell>0$.
For interior $X$ in this cylinder, $d(X+te)=d(X)$ whenever
$|t|\le|h|<\frac{\rho}{2}$. We integrate the Hessian bound along this segment and use $s^{k+1}\le Cd$ at its center. Uniformly for $0<|h|<\frac{\rho}{2}$,
\[
 \begin{gathered}
 0\le Q_h[u](X)\le C(1+d(X)^{-q})
 \le C(1+s(X)^{-m}),\\
 m=q(k+1).
 \end{gathered} \tag{B}
\]

Put $D=\rho^2-|z-z_*|^2$, $\phi=D^{-2m}$ and
$F=\phi(g_0-r-r^2)$, with $g_0-r-r^2\ge1$ on the cylinder. Set
$L_u=u^{ij}\partial_{ij}+\frac{k}{u_r}\partial_r$.
By (A),
\[
 |u^{r\alpha}\phi_\alpha|\le\sqrt{u^{rr}}\sqrt{u^{\alpha\beta}\phi_\alpha\phi_\beta}
 \le C\sqrt{u^{rr}}\,\frac{\phi}{D}.
\]
Using $|D_z^2\phi|\le \frac{C\phi}{D^2}$ and Young's inequality, we obtain,
with $s=u_r$, $L_uF\le\phi\bigl(\frac{C}{D^2}-\frac{k}{s}\bigr)$.
Choose $\sigma>0$ so that $L_uF<0$ for $0<s<\sigma D^2$ and
$\sigma\rho^4<c_\ell$. This region is disjoint from $\{r=r_T\}$. On $s=\sigma D^2$ we have $Q_h[u]\le CF$ by (B).

Fix $h\ne0$. As $D\downarrow0$, we have $F\to\infty$, while
$Q_h[u]\le\frac{2\|Du\|_\infty}{|h|}$.
On the graph, $\partial_rQ_h[u]=0$ and $F_r<0$. In the interior,
$L_uQ_h[u]\ge0$ by concavity. The maximum principle implies
$Q_h[u]\le CF$ throughout the region. At $z=z_*$, this bound is uniform in $h$.
For $s\ge\sigma\rho^4$, use (B). Finally, let $h\to0$.
\end{proof}

\section{Global comparison and quadratic splitting}\label{sec:splitting}\label{extR:conditional-splitting}
We compare the second derivatives of $f$ and $V$ to prove that
$D^2f\,e_y$ is constant. We first prove the splitting for one Euclidean
variable.

\begin{theorem}
\label{extR:split-conditional}
Let $f$ satisfy the hypotheses of Proposition~\ref{extR:globalpair},
and set $C_0=\{r>h(z)\}$. There are constants $N_0>0$,
$\beta\in\R^p$ and $\gamma\in\R$ such that
\[
 f(r,z,y)=F(r,z)+\frac{N_0}{2}(y+\beta\cdot z+\gamma)^2,
\]
where $F(r,z)=\inf_y f(r,z,y)\ge0$,
$F\in C^\infty(C_0)\cap C^1_{\loc}(\overline C_0)$, $D^2F>0$, and
\[
 \begin{gathered}
 F_r^k\det D^2F=\frac{K}{N_0},\\
 F_r|_{\partial C_0}=0.
 \end{gathered}
\]
\end{theorem}
\begin{proof}
As in Proposition~\ref{extR:globalpair}, we first replace $f$ by
$f-L_0$, where $L_0(X)=f(0)+Df(0)\cdot X$.
Write $N=p+2$. By Proposition~\ref{loc:directional}, for each $a>0$ we have
\begin{equation}\label{extR:approximation}
 \begin{gathered}
 (f_H)_{yy}\le C_a \quad\text{on }\{f_H<a\},\\
 (v_H)_{\eta\eta}\le C_a \quad\text{on }\{v_H<a\},\quad \eta\perp e_s,\ |\eta|=1.
 \end{gathered}
\end{equation}

\smallskip\noindent\emph{Normalization of $\rho_H$.}
Fix $X^\circ\in C$ and $P^\circ\in\mathcal H$. Throughout the proof, take
$H\ge H_0>\max\{1,f(X^\circ),V(P^\circ)\}$.
For the matrices in Lemma~\ref{extR:polarity}, put $t_H=\frac{\rho_H^2}{H}$.
With $a=1$ in \eqref{extR:approximation},
$t_Hf_{yy}\le C$ on $\{f<H\}\cap C$. At $P^\circ$,
\[
 \frac{\lambda_{\min}(D^2V(P^\circ))}{t_H}
 \le (v_H)_{qq}(\frac{A_H^TP^\circ}{H})\le C,
\]
because the $q$ component of $A_H^{-T}e_q$ is $\frac{1}{\rho_H}$.
At $X^\circ$, we have $t_Hf_{yy}(X^\circ)\le C$.
Since $f_{yy}(X^\circ)>0$, $c\le t_H\le C$.
We multiply the last column by $\frac{\sqrt H}{\rho_H}$ and obtain
\[
 A_H=\begin{pmatrix}d_H&b_H^T&0\\0&M_H&0\\0&c_H^T&\sqrt H\end{pmatrix}.
\]
The rescaling and its inverse are uniformly bounded and preserve
the directions $e_r,e_y$. The inverse transpose leaves $s$ unchanged.
Hence the normalized estimates hold, and $0<f_{yy}\le M$ in $C$.

\smallskip\noindent\emph{Smooth dual limits.}
For $H_j\to\infty$, set $f_j(X)=H_j^{-1}f(A_jX)$,
$v_j(P)=H_j^{-1}V(H_jA_j^{-T}P)$ and $B_j=\frac{b_j}{d_j}$.
By Lemmas~\ref{extR:polarity} and~\ref{audit:interior}, a diagonal subsequence satisfies
\[
 \begin{gathered}
 \kappa_j\to\kappa_\infty>0,\\
 v_j\to v_\infty\quad \hbox{in }C^m_{\loc}(\mathcal H)\quad (m<\infty).
 \end{gathered}
\]
Since $cI\le D^2v_j\le CI$ on compact subsets of $\mathcal H$,
$v_\infty$ is smooth and strictly convex, and
$\det D^2v_\infty=\frac{s^k}{\kappa_\infty}$.

\smallskip\noindent\emph{Comparison in $C$.}
Set $N_0=\sup_Cf_{yy}\in(0,\infty)$. For $\tau\ne0$, define
\[
 w_\tau(X)=\frac{f(X+\tau e_y)+f(X-\tau e_y)-2f(X)}{\tau^2}.
\]
Then $0\le w_\tau\le N_0$. Let $a_0=\frac{N+k}{k+2}$ and fix
$0<\delta<\frac{1}{a_0}$ independently of $h>0$. Set
$\Omega_h=\{rf_r<\delta(h+f)\}$ and
$\Sigma_h=\{rf_r=\delta(h+f)\}$.
We have $\Sigma_h\ne\varnothing$. Indeed, $g(r)=f(r,0,0)$ satisfies
$g(0)=g'(0)=0$, $g(r)\to\infty$ and $rg'(r)\ge g(r)$ by convexity.
Hence $rg'(r)-\delta(h+g(r))$ changes sign. Since $w_\tau$ is bounded,
we apply Lemma~\ref{loc:front-comparison} with $\ell=r$,
$m=\sup_{\Sigma_h}w_\tau$ and $\lambda=0$ to obtain
$w_\tau\le\sup_{\Sigma_h}w_\tau$ on $\Omega_h$.

Choose $Z_j\in C$ such that $f_{yy}(Z_j)>N_0-\frac{1}{j}$.
Since $w_\tau(Z_j)\to f_{yy}(Z_j)$ as $\tau\to0$, choose
$0<\tau_j<\frac{1}{j}$ with $w_{\tau_j}(Z_j)>N_0-\frac{2}{j}$.
Take $h_j\ge j$ so that $Z_j\in\Omega_{h_j}$. By comparison, there is
$X_j\in\Sigma_{h_j}$ such that $w_{\tau_j}(X_j)>N_0-\frac{3}{j}$.
Set $H_j=h_j+f(X_j)$. Then $f(X_j)\le H_j$ and
$r_jf_r(X_j)=\delta H_j$.

Set $Y_j=A_j^{-1}X_j$ and $P_j=Df_j(Y_j)$. By superlinear growth,
convexity and the local bounds, $|Y_j|+|P_j|\le C$.
Moreover, $[(Y_j)_r+B_j\cdot(Y_j)_z](P_j)_s=\delta$.
Hence $(P_j)_s\ge c>0$, and $P_j\to P_\infty\in\mathcal H$ along a subsequence.

By $Y_j=Dv_j(P_j)$ and Lemma~\ref{audit:interior}, there is $d>0$,
independent of $j$, such that $B_d(Y_j)\subset A_j^{-1}C$ and
$D^2f_j\le CI$ on $B_d(Y_j)$.
Set $\widehat\tau_j=\frac{\tau_j}{\sqrt{H_j}}\to0$. Since
$A_je_y=\sqrt{H_j}e_y$, the mean value theorem implies that, for large $j$,
there is $\widetilde Y_j\in Y_j+[-\widehat\tau_j,\widehat\tau_j]e_y$ with
$w_{\tau_j}(X_j)=(f_j)_{yy}(\widetilde Y_j)$.
For $\widetilde P_j=Df_j(\widetilde Y_j)$, we have
$|\widetilde P_j-P_j|\le C\widehat\tau_j\to0$.

Define $\chi_j(P)=e_q^T(D^2v_j(P))^{-1}e_q$ and
$\chi(P)=e_q^T(D^2v_\infty(P))^{-1}e_q$.
By the identity for the inverse Hessian, $0<\chi_j\le N_0$ on
$\mathcal H$ and $N_0-\frac{3}{j}<\chi_j(\widetilde P_j)\le N_0$.
By smooth convergence near $P_\infty$ and
$\widetilde P_j\to P_\infty\in\mathcal H$, we obtain
$\chi(P_\infty)=N_0$ and $\chi\le N_0$ in $\mathcal H$.

Write $M=D^2v_\infty$, $J=M^{-1}$, $a=Je_q$ and
$T_{ij}=\sum_l a_l(v_\infty)_{ijl}$. Differentiating $JM=I$ and $\log\det M=k\log s-\log\kappa_\infty$ twice, we find
\[
 L_{v_\infty}\chi
 =\operatorname{tr}(JTJT)-D^2(\log\det M)[a,a]
 =\operatorname{tr}(JTJT)+\frac{k a_s^2}{s^2}\ge0.
\]
By the strong maximum principle
\cite[Theorem~3.5]{GilbargTrudinger2001}, we obtain
$\chi\equiv N_0$ on the connected half space.
Both nonnegative terms vanish, so $T=0$ and $a_s=0$.

By symmetry of third derivatives,
\[
 \partial_i a=-J(\partial_i M)a=-JT e_i=0.
\]
Thus $a$ is constant. Set $\eta'=\frac{a}{N_0}$. This vector is horizontal
and has $q$ component one. Since $M\eta'=\frac{e_q}{N_0}$, we have
$(v_\infty)_{\eta'\eta'}=N_0^{-1}$ in $\mathcal H$.

\smallskip\noindent\emph{Comparison for the original potential.}
Set $\eta_j=\sqrt{H_j}A_j^{-T}\eta'$. These directions are horizontal
and have $q$ component one. At $P^\circ$, by \eqref{extR:approximation},
\[
 V_{\eta_j\eta_j}(P^\circ)
 =(v_j)_{\eta'\eta'}(\frac{A_j^TP^\circ}{H_j})\le C.
\]
Since $D^2V(P^\circ)>0$, the directions are bounded.
We fix a subsequence such that $\eta_j\to\eta_\infty$.
With $c_0=\frac{N}{k+2}<a_0$, define in $\mathcal H$
$G_j=H_j+V-c_0sV_s$, $O_j=\{sV_s<\delta(H_j+V)\}$ and
$\Gamma_j=\{sV_s=\delta(H_j+V)\}$.
Then $G_j\ge(1-c_0\delta)(H_j+V)$ on $\overline O_j$.
Fix $0<|\tau|\le1$ and define
\[
 Q_{j,\tau}(P)=\frac{V(P+\tau\eta_j)+V(P-\tau\eta_j)-2V(P)}{\tau^2}.
\]
For $P\in\mathcal H$, set $H=H_j$, $R=2(H+V(P))$,
$P_R=\frac{A_R^TP}{R}$ and $w_R=\frac{A_R^T\eta_j}{R}$.
By the ellipsoid inclusions in Lemma~\ref{extR:polarity},
\[
 |w_R|=\frac{\sqrt H}{R}|(A_j^{-1}A_R)^T\eta'|
 \le C H^{-\frac{1}{2}}\Bigl(\frac{R}{H}\Bigr)^{-\beta}\le C.
\]
Since $v_R(P_R)<1$, we have $|P_R|\le C$ by superlinear growth.
Also $A_Re_r=d_Re_r$ implies $w_R\perp e_s$. Hence
$P_R+t\tau w_R\in\mathcal H\cap B_C\quad(|t|\le1)$.
The local bound for $v_R$ and \eqref{extR:approximation} imply
\[
 \begin{aligned}
 Q_{j,\tau}(P)
 &=R\int_{-1}^1(1-|t|)
 D^2v_R(P_R+t\tau w_R)[w_R,w_R]\,dt\\
 &\le CR|w_R|^2
 \le C\Bigl(1+\frac{V(P)}{H_j}\Bigr)^{1-2\beta}.
 \end{aligned}
\]
All constants here are independent of $j,P,\tau$.
For $s=0$, let $h\downarrow0$ at $P+h e_s$ with $\tau$ fixed and use continuity.

Let $\epsilon,\xi>0$. Since $1-2\beta<1$, there is $T>0$,
independent of $j,\tau$, such that
$Q_{j,\tau}\le\frac{\epsilon G_j}{2H_j}$ on
$\overline O_j\cap\{V\ge TH_j\}$.
On $\Gamma_j\cap\{V<TH_j\}$, set $P'=\frac{A_j^TP}{H_j}$.
Since $v_j(P')\le T$, we have $|P'|+|Dv_j(P')|\le C$.
On $\Gamma_j$,
\[
 P'_s\bigl[(v_j)_s+B_j\cdot D_\zeta v_j\bigr](P')=\delta(1+v_j(P'))\ge\delta.
\]
Thus $P'_s\ge c>0$, and $P'$ lies in a compact subset of
$\mathcal H$ independent of $j$. Since
$\frac{A_j^T\eta_j}{H_j}=\frac{\eta'}{\sqrt{H_j}}$, we have
\[
 Q_{j,\tau}(P)=\int_{-1}^1(1-|t|)
 (v_j)_{\eta'\eta'}(P'+\frac{t\tau\eta'}{\sqrt{H_j}})\,dt .
\]
By smooth convergence on a larger compact subset of $\mathcal H$,
for all large $j$ we have $Q_{j,\tau}\le N_0^{-1}+\xi$ on
$\Gamma_j\cap\{V<TH_j\}$, uniformly for $0<|\tau|\le1$.

Combining these estimates, we obtain
$Q_{j,\tau}\le N_0^{-1}+2\xi+\frac{\epsilon G_j}{H_j}$ on $\Gamma_j$.
For fixed $j,\tau$, the growth bound is $o(H_j+V)$ at infinity.
Apply Lemma~\ref{loc:front-comparison} with $\ell=r$, $h=H_j$,
$m=N_0^{-1}+2\xi$ and $\lambda=\frac{\epsilon}{H_j}$. Then
$Q_{j,\tau}\le N_0^{-1}+2\xi+\frac{\epsilon G_j}{H_j}$ in $O_j$.
Fix $P\in\mathcal H$ and $0<|\tau|\le1$. For all sufficiently large $j$,
$P\in O_j$, and $\frac{G_j(P)}{H_j}\to1$. Since $\eta_j\to\eta_\infty$ and
$V$ is continuous, letting $j\to\infty$ with $\tau$ fixed, we obtain
\[
 \frac{V(P+\tau\eta_\infty)+V(P-\tau\eta_\infty)-2V(P)}{\tau^2}
 \le N_0^{-1}+2\xi+\epsilon.
\]
First let $\epsilon,\xi\downarrow0$, then $\tau\to0$. Thus $V_{\eta_\infty\eta_\infty}\le N_0^{-1}$.

\smallskip\noindent\emph{Equality for the original potential.}
For $P=Df(X)$, set $A=D^2f(X)$ and
$\eta=\eta_\infty$. Since $\eta\cdot e_y=1$,
\[
 0\le (Ae_y-N_0\eta)^TA^{-1}(Ae_y-N_0\eta)
 =f_{yy}-2N_0+N_0^2V_{\eta\eta}\le0.
\]
Thus $D^2f\,e_y=N_0\eta$ in $C$. We now add $L_0$ and use $f$
for the original potential; the Hessian identity remains valid.
Write $\eta=(0,\beta,1)$. Since $C$ is connected, there is a constant
$\gamma$ such that $f_y=N_0(y+\beta\cdot z+\gamma)$.
Set $m(z)=-\beta\cdot z-\gamma$ and define
$F(r,z)=f(r,z,m(z))$, $(r,z)\in\overline C_0$.
Integrating from $m(z)$ to $y$ proves the representation in $C$.
By continuity, it holds on $\overline C$. Thus
$F=\min_y f=\inf_y f\ge0$ on $\overline C_0$, with unique minimizer
$y=m(z)$.

Since $m$ is affine, $F\in C^\infty(C_0)$. Also $f_y(r,z,m(z))=0$, so
$DF(r,z)=(f_r,D_zf)(r,z,m(z))$ extends continuously to $\overline C_0$.
For $\xi=(\xi_r,\xi_z)\ne0$, put
$L\xi=(\xi_r,\xi_z,-\beta\cdot\xi_z)$. Then
$D^2F[\xi,\xi]=D^2f[L\xi,L\xi]>0$ and $\det D^2f=N_0\det D^2F$.
The first identity holds on the minimizing graph. For the second,
we use the change of variables $y\mapsto y+\beta\cdot z+\gamma$,
whose Jacobian is one.
Since $m$ is independent of $r$, we have $F_r=f_r>0$ on the minimizing graph. Hence
$F_r^k\det D^2F=\frac{K}{N_0}$ in $C_0$.
For $(r,z)\in\partial C_0$, $(r,z,m(z))\in\partial C$, and
$F_r(r,z)=f_r(r,z,m(z))=0$.

\end{proof}

\begin{proof}[Completion of the proof of Theorem~\ref{rf:splitting}]
For $x\in(0,\infty)^l$, write $x=r\ee+z$ with $z\cdot\ee=0$. Then
$r>\max_i(-z_i)\ge0$ and $\partial_r=\sum_i\partial_{x_i}$.
In an orthonormal transverse basis, the squared Jacobian is $l$.
Apply Theorem~\ref{extR:split-conditional} to each $y$ coordinate of
$\Phi$, taking the other $y$ variables as graph coordinates.
The corresponding Hessian columns are constant and orthogonal to
$(\ee,0)$. Integrating, we obtain $D_y\Phi=Ay+Cx+d$,
where $A=A^T>0$, $C$ and $d$ are constant, with $C\ee=0$.
Set $B=A^{-1}C$ and $b=A^{-1}d$. Then $B\ee=0$ and
\[
 D_y\!\Bigl[\Phi-\frac{1}{2}(y+Bx+b)^TA(y+Bx+b)\Bigr]=0.
\]
This proves \eqref{intro:split-equation}.

Its affine minimizing graph is $y=-Bx-b$, so the residual is
$\psi=\inf_y\Phi\ge0$ and has the stated interior smoothness and
$C^1$ regularity up to the boundary. By the Schur complement formula, $D^2\psi>0$.
In the original coordinates the affine shear
$(x,y)\mapsto(x,y+Bx+b)$ has unit Jacobian, and $D^2\Phi$ is the
pullback of $\operatorname{diag}(D^2\psi,A)$. Hence
$\det D^2\Phi=(\det A)\det D^2\psi$ and $W=\sum_i\psi_{x_i}$.
Substitution in \eqref{intro:real-equation} proves
\eqref{intro:residual-equation}, including the boundary trace of the weight.
Conversely, from the same identities we verify that any stated choice of
$\psi,A,B,b$ defines a strictly convex solution of
\eqref{intro:real-equation}. The function $\Phi$ has the required regularity and satisfies $\Phi\ge\psi\ge0$.
\end{proof}

\clearpage
\section*{Acknowledgments}
This work was supported by the National Key R\&D Program of China (2025YFA1017603).

\section*{Declaration of AI use}
During preparation the authors used ChatGPT (OpenAI) for language editing,
exploration of proof strategies and calculations.

\begingroup\footnotesize
\par\medskip\noindent
Jiahuan Li, Yilu Liu and Xi-Nan Ma. School of Mathematical Sciences,
University of Science and Technology of China, Hefei, Anhui 230026,
P.R. China.\\
\href{mailto:jiahuan@mail.ustc.edu.cn}{jiahuan@mail.ustc.edu.cn},
\href{mailto:liuylgeoanaly@mail.ustc.edu.cn}{liuylgeoanaly@mail.ustc.edu.cn},
\href{mailto:xinan@ustc.edu.cn}{xinan@ustc.edu.cn}, respectively.
\par\smallskip\noindent
Dekai Zhang. School of Mathematical Sciences, Key Laboratory of Mathematics
and Engineering Applications (Ministry of Education), Shanghai Key Laboratory
of PMMP, East China Normal University, Shanghai 200241, China.
\href{mailto:dkzhang@math.ecnu.edu.cn}{dkzhang@math.ecnu.edu.cn}.
\par\endgroup

\end{document}